\documentclass[a4paper]{article}
\usepackage{amsmath,mathtools,amsthm,a4wide}
\usepackage{amssymb, authblk}
\usepackage[mathscr]{eucal}
\usepackage{bm,booktabs}
\usepackage{bbm}
\usepackage[shortlabels]{enumitem}
\usepackage{hyperref}

\usepackage{amsmath}
\usepackage{graphicx,tabularx}
\allowdisplaybreaks[1]

\newcommand{\mc}{\multicolumn}
\newcommand{\pr}{\rightarrow}

\newcommand{\ba}{\begin{array}}
\newcommand{\ea}{\end{array}}
\newcommand{\be}{\begin{equation}}
\newcommand{\ee}{\end{equation}}

\newcommand{\vart}{\vartheta}
\newcommand{\varp}{\varphi}
\newcommand{\eps}{\varepsilon}

\newcommand{\il}{\int\limits}
{\begin{list}{}{\setlength{\rightmargin}{0cm}
                \setlength{\listparindent}{0cm}
                \settowidth{\labelwidth}{\mbox{#1}}
                \setlength{\leftmargin}{1.1\labelwidth}
                \setlength{\labelsep}{.1\labelwidth}}}%
{\end{list}}

\newcommand{\ei}{\end{inspring}}

\newcommand{\bfx}{{\bf x}}

\newcommand{\beq}{\begin{equation}}
\newcommand{\eq}{\end{equation}}
\catcode`\@=11

\font\tenmsa=msam10 \font\sevenmsa=msam7 \font\fivemsa=msam5
\font\tenmsb=msbm10 \font\sevenmsb=msbm7 \font\fivemsb=msbm5
\newfam\msafam
\newfam\msbfam
\textfont\msafam=\tenmsa  \scriptfont\msafam=\sevenmsa
  \scriptscriptfont\msafam=\fivemsa
\textfont\msbfam=\tenmsb  \scriptfont\msbfam=\sevenmsb
  \scriptscriptfont\msbfam=\fivemsb

\def\Bbb{\ifmmode\let\next\Bbb@\else
 \def\next{\errmessage{Use \string\Bbb\space only in math mode}}\fi\next}
\def\Bbb@#1{{\Bbb@@{#1}}}
\def\Bbb@@#1{\fam\msbfam#1}
\newcommand{\dR}{{\Bbb R}}

\newcommand{\dC}{{\Bbb C}}

\newcommand{\dE}{{\Bbb E}}

\newtheorem{thm}{Theorem}

\newtheorem{assumption}[thm]{Assumption}

\begin{document}
\title{Heavy-traffic single-server queues and\\ the transform method}
\author{M.A.A. Boon, A.J.E.M. Janssen and J.S.H. van Leeuwaarden}

\maketitle

\begin{abstract}
Heavy-traffic limit theory deals with queues that operate close to criticality and face severe queueing times. Let $W$ denote the steady-state waiting time in the ${\rm GI}/{\rm G}/1$ queue. Kingman (1961) showed that $W$, when appropriately scaled, converges in distribution to an exponential random variable as the system's load approaches 1. The original proof of this famous result uses the transform method. Starting from the Laplace transform of the pdf of $W$ (Pollaczek's contour integral representation), Kingman showed convergence of transforms and hence weak convergence of the involved random variables. We apply and extend this transform method to obtain convergence of moments with error assessment. We also demonstrate how the transform method can be applied to so-called nearly deterministic queues in a Kingman-type and a Gaussian heavy-traffic regime. We demonstrate numerically the accuracy of the various heavy-traffic approximations.

%\noindent\textbf{Keywords: }G/G/1 queue, heavy-traffic analysis, approximations, error assessment
\end{abstract}

\maketitle

%\tableofcontents

\section{Introduction and results} \label{sec1}

{The title of this contribution to the memorial issue for J.W.~Cohen refers to the The Single-Server Queue, the monumental book \cite{cohen2012single} in which J.W.~Cohen teaches the reader how to use complex analysis and transform methods to obtain mathematical rigorous results for the general GI/G/1 queue and its many extensions. In turn, J.W.~Cohen admired the work of F.~Pollaczek, in particular for the analytic treatment of queues by means of complex function theory and integral equations \cite{cohen1993complex}, techniques that also  feature prominently in The Single-Server Queue, and in this paper on the GI/G/1 queue in heavy traffic. }

F.~Pollaczek initiated the analysis of the GI/G/1 queue
%and GI/G/$s$ queue
in the 1940s and 1950s, and obtained a contour integral representation for the Laplace transform of the steady-state waiting time.  J.F.C.~Kingman introduced heavy-traffic analysis in the 1960s \cite{ref10,ref11}.
For the {\rm GI}/{\rm G}/1 queue in a regime where the system load tends to 1, Kingman showed that, under certain conditions, the Laplace transform of the pdf of an appropriately normalized steady-state waiting time (Pollaczek's contour integral representation) converges to the Laplace transform of an exponential distribution. We will refer to this technique---to show convergence in distribution by convergence of transforms---as the transform method.
To explain Kingman's result in more detail, let $W$
denote the steady-state waiting time in the {\rm GI}/{\rm G}/1 queue, which solves the stochastic equation
\beq \label{1.1}
W\stackrel{d}{=} \Big(W+V-\frac{1}{\rho}\,U\Big)^+,
\eq
with $x^+=\max\{0,x\}$. Here, $V$ is the generic service time with mean 1 and variance $\sigma_V^2\in(0,\infty)$, $U$ is the generic inter-arrival time with mean 1 and variance $\sigma_U^2\in(0,\infty)$ and $\rho\in(0,1)$. It is assumed that $V$ and $U$ are independent. Since convergence of transforms implies convergence of distributions, Kingman effectively showed using the transform method that, for $W=W_{\alpha}$ solving (\ref{1.1}) with $\alpha=1-\rho$,
\beq \label{1.2}
P(\alpha\,W_{\alpha}\leq t)\pr 1-e^{-2t/\sigma^2},~~~~~~\alpha\downarrow 0,
\eq
for all $t\geq0$ with $\sigma^2=\sigma_U^2+\sigma_V^2$.

Kingman's observation that heavily loaded systems admit a simple scaling limit triggered a surge of research in the 1960s and 1970s; see \cite{ref3,ref6,ref4,ref15,ref19,ref12,ref13}, among others. In the decades that followed, heavy-traffic analysis, and more generally, stochastic-process limits, developed into popular topics in the applied probability community, with queueing theory as one of its many applications. The general idea remained to consider a parametrized set of systems and to identify the limiting system as the parameter converges to a limiting value yielding criticality (e.g.\ $\alpha\downarrow0$ as in (\ref{1.2})).

Kingman's transform method is thus largely based on Pollaczek's contour integral that we introduce next, see \cite{ref1}. Assume analyticity of $\psi(s)=\dE\,[\exp(s(V-\frac{1}{\rho}\,U))]$ in an open strip containing $|\mbox{Re}(s)|\leq\delta$ for some $\delta>0$; in particular, all moments of $V-\frac{1}{\rho}\,U$ are finite. Then, Pollaczek's integral reads
\beq \label{1.3}
\dE\,[e^{-sW}]=\exp\,\Bigl\{\frac{-1}{2\pi i}\,\il_C\,\frac{s}{z(s-z)}\,\log(1-\psi({-}z))dz\Bigr\},
\eq
where $C$ is a contour to the left of, and parallel to, the imaginary axis, and to the right of the singularities of $\log(1-\psi({-}z))$ in the left half-plane, and $s$ is any complex number to the right of $C$. Kingman uses (\ref{1.3}) to prove that $\dE\,[\exp({-}s\alpha W)]\pr(1+\sigma^2s/2)^{-1}$ in $\mbox{Re}(s)\geq0$ as $\alpha\downarrow 0$, yielding (\ref{1.2}).

Being tailor-made for the steady-state ${\rm GI}/{\rm G}/1$ queue, the transform method that uses contour integral representations has rarely been applied to other queueing models. A notable exception is the heavy-traffic analysis of O.J.~Boxma and J.W.~Cohen \cite{ref3b} for the ${\rm GI}/{\rm G}/1$ queue with heavy-tailed distributions, so when the second moment of the service time and/or interarrival time is infinite. Boxma and Cohen apply the transform method to an extended form of Pollaczek's integral \eqref{1.3} to identify the heavy-traffic limit.
Other studies that apply this transform method for heavy-traffic analysis are  \cite{ref12,ref13} on the GI/G/$s$ queue and \cite{boon2019pollaczek,boon2021optimal} on the fixed-cycle traffic-light queue, a variation of the GI/G/1 queue.

More probabilistic methods for proving heavy-traffic results developed later use functional limit theorems, and typically establish that the sample path of a scaled waiting time converges to some limiting stochastic process. One is then faced with the problem of showing that the steady-state of the limiting process corresponds to the limiting steady-state of the queueing model in heavy traffic. This requires and interchange-of limits argument, which  is often challenging as it involves proving tightness of sequences of probability measures. The transform method works directly with the steady-state random variables, and thus avoids the problem of interchanging limits.

\subsection{Nearly deterministic queues and the transform method} \label{subsec1.1}

Next to the classical heavy-traffic setting, we will consider so-called nearly deterministic ${\rm GI}_n/{\rm G}_n/1$ queues, whose heavy-traffic behavior has been investigated in \cite{ref17,ref18} using sample-path methods. Nearly-deterministic queues are motivated by cycling thinning. To explain this, denote for all $n=1,2,...$
\beq \label{1.4}
W_n\stackrel{d}{=} \Bigl(W_n+V_n-\frac{1}{\rho_n}\,U_n\Bigr)^+,
\eq
with
\beq \label{1.5}
V_n=\frac1n\,\sum_{k=1}^n\,V_{n,k},~~~~~~U_n=\frac1n\,\sum_{k=1}^n\,U_{n,k},
\eq
where $V_{n,k}$ are i.i.d.\ copies of $V$ and $U_{n,k}$ are i.i.d.\ copies of $U$, with $V$ and $U$ as before, and $\rho_n\in(0,1)$.
The cyclic thinning thus regards each interarrival (service) time as the $n$-th occurrence in a sequence of i.i.d.~random variables, which mitigates the random fluctuations. For instance, if this sequence consists of exponential random variables, the interarrival or service time would follow an Erlang distribution, and if $n$ is large, this becomes more and more `deterministic'.

Interesting heavy-traffic regimes now arise when $\rho_n\pr 1$ as $n\pr\infty$. In \cite{ref17,ref18}, two heavy-traffic regimes are considered. The first (Kingman-type) regime assumes that $(1-\rho_n)n\pr\beta$ as $n\pr\infty$ for some fixed $\beta>0$. In this case, $W_n$ converges in distribution to an exponential random variable with mean $\sigma^2/2\beta$, where $\sigma^2=\sigma_U^2+\sigma_V^2$. The second (Gaussian) regime assumes $(1-\rho_n)\,\sqrt{n}\pr\beta$ as $n\pr\infty$ for some fixed $\beta>0$. In this case $\sqrt{n}\,W_n/\sigma$ converges in distribution to the all-time maximum $M_{\beta}$ of a one-dimensional directed random walk, starting at 0, with normally distributed step sizes of mean $-\beta$. This Gaussian random walk and $M_{\beta}$ are well studied, see \cite{ref2,ref5,ref7,ref8,ref16}.

\subsection{Main results} \label{subsec1.2}

In the present paper we apply the transform method for heavy-traffic analysis of the ${\rm GI}/{\rm G}/1$ and ${\rm GI}_n/{\rm G}_n/1$ queues.
%Let us now summarize the main theoretical results; we discuss the numerical implications in Subsection \ref{subsec1.num}.
We first  consider the classical heavy-traffic regime, and provide a detailed proof of a version of Kingman's original result using the transform method. We do this with  service times $V$ and interarrival times $U$ that do not depend on $\rho=1-\alpha$ (in Kingman's original result, a controlled dependence of $V$ and $U$ on $\alpha$ is allowed). In this more restricted setting, we show the following.

\begin{thm} \label{thm1}
With $\sigma_{\alpha}^2=(\sigma_V^2+\rho^{-2}\sigma_U^2)\rho$,
\beq \label{1.6}
\dE\,[e^{-\alpha sW}]=(1+\sigma_{\alpha}^2\,s/2)^{-1}+O(\alpha\,\log(1/\alpha)),~~~~~~\alpha\downarrow0,
\eq
uniformly in any bounded set of $s$ with ${\rm Re}(s)\geq{-}1/2\sigma_{\alpha}^2$.
\end{thm}
As a consequence of Theorem~\ref{thm1}, we have for any $k=1,2,...$
\beq \label{1.7}
\dE\,[(\alpha W)^k]=k!(\tfrac12\sigma_{\alpha}^2)^k+O(\alpha\,{\log}(1/\alpha)),~~~~~~\alpha\downarrow 0.
\eq
We observe that for $k=1$ in \eqref{1.7} we get
\beq \label{1.7a}
\dE\,[\alpha W]=\tfrac12\sigma_{\alpha}^2+O(\alpha\,{\log}(1/\alpha)),~~~~~~\alpha\downarrow 0.
\eq
It turns out that, after appropriate identifications, the quantity $\tfrac12\sigma_\alpha^2$ at the right-hand side of \eqref{1.7a} coincides with the right-hand side of (6) in \cite{ref5a} (Kingman's bound, see \cite{ref11a}, Theorem~2 and (33)). In \cite{ref5a} there are discussed sharpenings of Kingman's bound, including a new upper bound \cite[Equation (17)]{ref5a}, involving Lambert's W function.

We next apply the transform method to the ${\rm GI}_n/{\rm G}_n/1$ queues described in Subsection~\ref{subsec1.1}, where we let $(1-\rho_n)\asymp 1/n$, meaning that there are fixed $\beta_1$ and $\beta_2$ with $0<\beta_1\leq\beta_2<\infty$ such that $(1-\rho_n)n\in[\beta_1,\beta_2]$, $n=1,2,...\,$. This leads to the following result.
\begin{thm} \label{thm2}
With $\gamma_n=(\sigma_V^2+\rho_n^{-2}\sigma_U^2)\rho_n/(2n(1-\rho_n))$,
\beq \label{1.8}
\dE\,[e^{-tW_n}]=(1+\gamma_nt)^{-1}+O\Bigl(\frac{{\rm log}\,n}{\sqrt{n}}\Bigr),~~~~~~n\pr\infty,
\eq
uniformly in any bounded set of $t$ with ${\rm Re}(t)\geq{-}1/4\gamma_n$.
\end{thm}
From Theorem~\ref{thm2} we obtain for any $k=1,2,...$
\beq \label{1.9}
\dE\,[W_n^k]=k!\,\gamma_n^k+O\Bigl(\frac{{\rm log}\,n}{\sqrt{n}}\Bigr),~~~~~~n\pr\infty.
\eq

We then proceed to apply the transform method to the ${\rm GI}_n/{\rm G}_n/1$ when we let $(1-\rho_n)\asymp 1/\sqrt{n}$, meaning that there are fixed $\beta_1$ and $\beta_2$ with $0<\beta_1\leq\beta_2<\infty$ such that $(1-\rho_n)\,\sqrt{n}\in[\beta_1,\beta_2]$, $n=1,2,...\,$. In this regime, the integration contour $C$ occurring in (\ref{1.3}) can be chosen to pass through the saddle point $z=\zeta_{sp}$ of $h(z)={\rm log}(\psi({-}z))$ on the negative real axis, allowing a saddle point analysis  (under an additional assumption). We show the following, with $M_{\beta}$ as in Subsection~\ref{subsec1.1}.
\begin{thm} \label{thm3}
With $\sigma_n=(h''(\zeta_{sp}))^{1/2}$ and $\beta_n={-}\zeta_{sp}\sigma_n\sqrt{n}$,  $\sigma_n\asymp1$, $\beta_n\asymp 1$ as $n\pr\infty$, and
\beq \label{1.10}
\dE\,\Bigl[\exp\Bigl({-}s\,\frac{\sqrt{n}}{\sigma_n}\,W_n\Bigr)\Bigr]=\dE\,[\exp({-}sM_{\beta_n})]+O\Bigl(\frac{1}{\sqrt{n}}\Bigr),~~~~~~n\pr\infty,
\eq
uniformly in any bounded set of $s$ with ${\rm Re}(s)\geq{-}\,\frac12\,\beta_n$.
\end{thm}
From Theorem~\ref{thm3} we get for any $k=1,2,...$
\beq \label{1.11}
\dE\,\Bigl[\Bigl(\frac{\sqrt{n}}{\sigma_n}\,W_n\Bigr)^k\Bigr]=m_k(\beta_n)+O\Bigl(\frac{1}{\sqrt{n}}\Bigr),~~~~~~n\pr\infty,
\eq
where $m_k(\beta)=\dE\,[M_{\beta}^k]$. Theorem~\ref{thm3} can be refined by taking account of $h'''(\zeta_{sp})$ in the saddle point analysis. This yields Theorem~\ref{thm4}, see Section~\ref{sec5} for its precise formulation, where the right-hand side of (\ref{1.10}) is replaced by $\dE\,[\exp({-}R_n(s)\,M_{B_n})]+O(1/n)$, with appropriate non-linear transforms $R_n(s)$ and $B_n(\beta)=B_n$ of $s$ and $\beta$. Theorem~\ref{thm3} and its refinement Theorem~\ref{thm4} result from an adaptation of the saddle point method developed in \cite{ref9}, \cite{ref14}.
Theorem \ref{thm5} in Section \ref{sec5} gives a consequence of Theorem~\ref{thm4} on the level of moments.
%It is worth mentioning that the $O(1/\sqrt{n})$ at the right-hand sides of \eqref{1.10} and \eqref{1.11} can be replaced by $O(1/n)$ when the third cumulants of $V$ and $U$ are equal.

Theorems~\ref{thm1}-\ref{thm4} generalize and refine some classical heavy-traffic results. Theorem~\ref{thm1} recovers Kingman’s weak convergence result \eqref{1.2}, and generalizes
this to a heavy-traffic limit theorem for all moments of $W$. The refined heavy-traffic approximations not only provide the heavy-traffic limit, but also contain pre-limit information (for $\rho$ away from 1) and shed light on the rate of convergence (the speed at which the limit is attained, as a function of the scaling parameter). Similarly, Theorems~\ref{thm2}-\ref{thm4}
recover results in \cite{ref17}, \cite{ref18} for convergence in
distribution, and convergence of the first two moments. In the paper we show that
all normalized moments of $W_n$ converge to those of $M_\beta$, again with rate of convergence and refinements. As a consequence, the theoretical results in  Theorems~\ref{thm1}-\ref{thm4}  give sharp approximations, not only in heavy traffic, but also in more moderate conditions.

We  demonstrate the accuracy of the approximations by comparing with exact results. We also address the computational aspects of numerically calculating complex contour integrals, which is required for both the exact and approximate performance analysis.
In particular, we explain how to calculate reliably the Pollaczek contour integrals with numerical integration. Since we operate in heavy-traffic regimes, numerical integration can become cumbersome, with integration contours closely passing the origin, %and/or integrands varying violently,
but we show how to deal with this.

\subsection{Organization of the paper} \label{subsec1.3}

In Section~\ref{sec2} we present assumptions and preliminaries on the function $\psi({-}z)=\dE\,[\exp({-}z(V-\frac{1}{\rho}\,U))]$ that occus in the various Pollaczek integrals in Sections~\ref{sec3}, \ref{sec4} and \ref{sec5}. Section~\ref{sec2} also contains information on the function $h(z)={\rm log}\,\psi({-}z)$ that is heavily used in the dedicated saddle point method of Section~\ref{sec5}; we avoid using saddle points in Sections~\ref{sec3}, \ref{sec4} on the Kingman-type results. In Section~\ref{sec3} we present the formulation and proof of our version of Kingman's classical result (Theorem~\ref{thm1}), yielding convergence (with error assessment) of the mgf and all moments of those of an exponentially distributed random variable with a tailored $\alpha$-dependent mean. In Section~\ref{sec4}, we consider the nearly deterministic queue in the regime $(1-\rho)\asymp 1/n$, and prove that the mgf and all moments of $W_n$ converge to those of a specifically designed exponentially distributed random variable (Theorem~\ref{thm2}). In Section~\ref{sec5} we present Theorem~\ref{thm3} and its refinement Theorem~\ref{thm4}, with consequences for moment convergence, for the nearly deterministic queue when  $(1-\rho)\asymp 1/\sqrt{n}$. This requires an additional assumption on the function $\psi$, allowing a saddle point approach to Pollaczek's integral around the saddle point $\zeta_{sp}$, at the expense of an exponentially small error. The proof of the refinement Theorem~\ref{thm4}, and its consequence (\ref{5.6}) %, (\ref{5.7})
for moment convergence, is rather involved and technical, so that we have deferred details to Appendices~A and B. In Section~\ref{sec6} we summarize in detail the computational schemes for the quantities we want to calculate via Pollaczek's formula \eqref{1.3}.
In Section~\ref{sec8} we present our conclusions.
In the appendix we illustrate the results of Sections~\ref{sec3}, \ref{sec4} and \ref{sec5} for the several specific cases, including $V$ and $U$ being Gamma distributed, and compare our asymptotic results with the results of numerical integration.

\section{Preliminaries} \label{sec2}

The convergence results of the Kingman type given in the present paper will be shown under the condition that the function $\psi$ is analytic in an open set containing a strip $|\mbox{Re}(z)|\leq \delta$ with $\delta>0$. For the convergence results for nearly deterministic queues related to the maximum of the Gaussian random walk, we require one more condition, additional to the analyticity assumption on $\psi$. The latter results are shown using the dedicated saddle point method, and this requires an assumption that guarantees that one can confine attention to the immediate vicinity of the saddle point on the negative real axis when conducting an asymptotic analysis on the relevant Pollaczek integral when $n\to\infty$. We refer to Section \ref{sec5} for the technical details. This additional condition is satisfied for all special distributions listed above, except the deterministic and two-point distributions that yield functions $\psi(-x+iy)$, $y\in\dR$, that are (nearly) periodic in $y$, and thus take values (near to) $\psi(-x)$ for certain $y$ away from 0.

We use in the sequel the letter $\delta$ to denote a generic positive number that may take case-dependent values. We consider random variables $V,U\geq0$ that are independent with $\dE\,[V]=1=\dE\,[U]$ and $0<\sigma_V^2+\sigma_U^2<\infty$. We furthermore assume that there is a $\delta>0$ such that $\dE\,[\exp(zV)]$, $\dE\,[\exp(zU)]$ are analytic in an open strip containing $-\delta\leq {\rm Re}(z)\leq\delta$. For $\rho\in(0,1)$, we let
\beq \label{2.1}
\psi({-}\zeta)=\psi({-}\zeta\,;\,\rho)=\dE\,[e^{-\zeta(V-\frac{1}{\rho}U)}].
\eq
Then $\psi({-}\zeta)$ is analytic in an open strip containing $-\delta\leq{\rm Re}(\zeta)\leq\delta$ for some $\delta>0$. Since
\beq \label{2.2}
\psi({-}\zeta)\equiv\il_{-\infty}^{\infty}\,e^{-\zeta t}\,d\,G(t),
\eq
with $G(t)$ the cumulative distribution function of $V-\frac{1}{\rho}\,U$, we have that $\psi({-}\zeta)$, $-\delta\leq\zeta\leq\delta$, is logarithmically convex. We have, uniformly in $\rho\in[\frac12,1]$,
\begin{eqnarray} \label{2.3}
\psi({-}\zeta) & = & \dE\,\Bigl[1-\zeta\Bigl(V-\frac{1}{\rho}\,U\Bigr)+\frac12\,\zeta^2\Bigl(V-\frac{1}{\rho}\,U\Bigr)^2\Bigr]+O(\zeta^3) \nonumber \\[3mm]
& = & 1+\Bigl(\frac{1}{\rho}-1\Bigr)\,\zeta+\frac12\,\Bigl(\sigma_V^2+\dfrac{1}{\rho^2}\,\sigma_U^2+\Bigl(1-\frac{1}{\rho}\Bigr)^2\Bigr)\,\zeta^2+O(\zeta^3), \quad |\zeta|\leq\delta,
\end{eqnarray}
for some $\delta>0$. Therefore, there is a $\delta>0$ such that $\psi({-}\zeta)\geq1/2$ when $\rho\in[\frac12,1]$ and $-\delta\leq\zeta\leq0$. Hence, by continuity, there is a $\delta>0$ such that for $\rho\in[\frac12,1]$
\beq \label{2.4}
h(\zeta)={\rm log}\,\psi({-}\zeta)={\rm log}(\dE\,[e^{-\zeta(V-\frac{1}{\rho}U)}])
\eq
is well-defined and analytic as a principal value logarithm in an open set containing the rectangle $-\delta\leq{\rm Re}(\zeta)\leq0$, $|{\rm Im}(\zeta)|\leq\delta$.

We have for $\rho\in[\frac12,1]$
\beq \label{2.5}
h(0)=0,~~~~~~h'(0)=\frac{1}{\rho}-1,~~~~~~h''(0)=\sigma_V^2+\frac{1}{\rho^2}\,\sigma_U^2,
\eq
and there is a $\delta>0$ such that
\beq \label{2.6}
h(\zeta)=\Bigl(\frac{1}{\rho}-1\Bigr)\,\zeta+\frac12\,\Bigl(\sigma_V^2+\frac{1}{\rho^2}\,\sigma_U^2\Bigr)\,\zeta^2+O(\zeta^3)
\eq
in an open set containing the rectangle $-\delta\leq{\rm Re}(\zeta)\leq0$, $|{\rm Im}(\zeta)|\leq\delta$. There is a $\delta>0$ such that the function $h(\zeta)$, $-\delta\leq\zeta\leq0$, is convex. Furthermore, $h$ has, for $\rho$ sufficiently close to 1, a unique saddle point $\zeta_{sp}\in[{-}\delta,0]$. We have
\beq \label{2.7}
\zeta_{sp}={-}\,\frac{1}{\rho}~\frac{1-\rho}{\sigma_V^2+\frac{1}{\rho^2}\,\sigma_U^2}+O((1-\rho)^2),
\eq
and
\beq \label{2.8}
h(\zeta_{sp})={-}\,\frac{1}{\rho^2}~\frac{(1-\rho)^2}{2(\sigma_V^2+\frac{1}{\rho^2}\,\sigma_U^2)}+O((1-\rho)^3),~~~~~h''(\zeta_{sp})=\sigma_V^2+\frac{1}{\rho^2}\,\sigma_U^2+O(1-\rho).
\eq

\section{Classical heavy-traffic result of the Kingman type} \label{sec3}

With $V$ and $U$ as in Section~\ref{sec2}, we let
\beq \label{3.1}
W\stackrel{d}{=}\Bigl(W+V-\frac{1}{\rho}\,U\Bigr)^+,
\eq
where $\rho=1-\alpha$ and $\alpha\downarrow 0$. We shall show that
\beq \label{3.2}
{\rm log}(\dE\,[e^{-\alpha sW}])={-}{\rm log}(1+\tfrac12\,\sigma_{\alpha}^2s)+O(\alpha\,{\rm log}(1/\alpha)),~~~~~~\alpha\downarrow0,
\eq
uniformly in any bounded set of $s$ with ${\rm Re}(s)\geq\frac12\,\mu_0$, where
\beq \label{3.3}
\sigma_{\alpha}^2=\frac{-1}{\mu_0}=\Bigl(\sigma_V^2+\frac{1}{\rho^2}\,\sigma_U^2\Bigr)\,\rho.
\eq
The mgf $\dE\,[\exp({-}t\bfx)]$ of an $\bfx$ having the exponential probability distribution $\vart\,e^{-\vart x}$, $x\geq0$, with mean $1/\vart$, is given by $(1+t/\vart)^{-1}$, ${\rm Re}(t)>{-}\vart$. Hence, all moments of $\alpha W$ are equal to the moments of the exponential probability distribution with mean $2/\sigma_{\alpha}^2$, with error $O(\alpha\,{\rm log}(1/\alpha))$ as $\alpha\downarrow 0$.

We show this result by using the Pollaczek result for $W$, in which we follow the argumentation as given by Kingman in \cite{ref10} rather closely. Observe that $V$ and $U$ are independent of $\alpha$, which allows us to be explicit about the error term in (\ref{3.2}).

From Pollaczek's result, we have
\beq \label{3.4}
{\rm log}(\dE\,[e^{-\alpha sW}])=\frac{-1}{2\pi i}\,\il_C\,\frac{\alpha s}{\zeta(\alpha s-\zeta)}\,{\rm log}(1-\psi({-}\zeta))d\zeta,
\eq
where $\psi({-}\zeta)$ is as in (\ref{2.1}), and $C$ is a line parallel to, and to the left, of the imaginary axis, and to the right of the singularities of ${\rm log}(1-\psi({-}\zeta))$ in the open left-half-plane, and $\alpha s$ in (\ref{3.4}) lies to the right of $C$. To be more detailed about the choice of $C$, we observe from (\ref{2.3}) that there is a $\delta>0$ such that
\beq \label{3.5}
\psi({-}\zeta)=1+\alpha\zeta/\rho+\tfrac12\,\sigma_{\alpha}^2\zeta^2/\rho+O(\alpha^2\zeta^2)+O(\zeta^3),~~~~~~|\zeta|\leq\delta,
\eq
where $\sigma_{\alpha}^2$ is as in (\ref{3.3}). The leading part $1+\alpha\zeta/\rho+\sigma_{\alpha}^2\zeta^2/2\rho$, considered for $\zeta\leq0$, in (\ref{3.5}) equals 1 for $\zeta=0$ and $\zeta=2\alpha\mu_0$, and is minimal for $\zeta=\alpha\mu_0$, where $\mu_0<0$ is as in (\ref{3.3}), with minimal value $1-\alpha^2/2\rho\sigma_{\alpha}^2$. By (\ref{2.2}) we have
\beq \label{3.6}
|\psi({-}\alpha\mu_0-i\eta)|\leq\psi({-}\alpha\mu_0)=1-\frac{\alpha^2}{2\rho\sigma_{\alpha}^2}+O(\alpha^3),~~~~~~\eta\in\dR,
\eq
and so
\beq \label{3.7}
|\psi({-}\alpha\mu_0-i\eta)|\leq1-\frac{\alpha^2}{4\rho\sigma_{\alpha}^2},~~~~~~\eta\in\dR,
\eq
when $\alpha$ is sufficiently small. For such $\alpha$, we can therefore choose $C=\{\alpha\mu_0+i\eta\;|\;\eta\in\dR\}$, with principal value of the log in (\ref{3.4}), and ${\rm Re|}(s)\geq\frac12\,\mu_0$.

In (\ref{3.4}) we substitute $\zeta=\alpha z$, with $z\in\{\mu_0+i\eta\;|\;\eta\in\dR\}=C_0$ and $d\zeta=\alpha\,dz$, to get
\beq \label{3.8}
{\rm log}(\dE\,[e^{-\alpha sW}])=\dfrac{-1}{2\pi i}\,\il_{C_0}\,\frac{s}{z(s-z)}\,{\rm log}(1-\psi({-}\alpha z))\,dz.
\eq
From (\ref{3.5}) we have when $|\alpha z|\leq\delta$,
\beq \label{3.9}
\psi({-}\alpha z)=1+\alpha^2z/\rho+\sigma_{\alpha}^2\alpha^2z^2/2\rho+O(\alpha^4z^4)+O(\alpha^3z^3),
\eq
and so
\beq \label{3.10}
\frac{1-\psi({-}\alpha z)}{-\alpha^2z/\rho}=1+\tfrac12\,\sigma_{\alpha}^2z+O(\alpha^2z)+O(\alpha z^2).
\eq
We have from (\ref{3.7}) that both $1-\psi({-}\alpha z)$ and $-\alpha z^2/\rho$ lie in the open right half-plane when $z\in C_0$ and $\alpha$ is sufficiently small. Hence, with principal-value logarithms,
\beq \label{3.11}
{\rm log}(1-\psi({-}\alpha z))={\rm log}\Bigl(\frac{1-\psi({-}\alpha z)}{-\alpha^2z/\rho}\Bigr)+{\rm log}({-}\alpha^2z/\rho),~~~~~~z\in C_0.
\eq
Since ${\rm Re}(s)\geq\frac12\,\mu_0$, we have by Cauchy's theorem
\beq \label{3.12}
\frac{-1}{2\pi i}\,\il_{C_0}\,\frac{s}{z(s-z)}\,{\rm log}({-}\alpha^2z/\rho)\,dz=0,
\eq
and so
\beq \label{3.13}
{\rm log}(\dE\,[e^{-\alpha sW}])=\frac{-1}{2\pi i}\,\il_{C_0}\,\frac{s}{z(s-z)}\,{\rm log}\Bigl(\frac{1-\psi({-}\alpha z)}{-\alpha^2z/\rho}\Bigr)\,dz.
\eq

We shall show now that
\beq \label{3.14}
{\rm log}\Bigl(\frac{1-\psi({-}\alpha z)}{-\alpha^2z/\rho}\Bigr)={\rm log}(1+\tfrac12\,\sigma_{\alpha}^2z)+O(\alpha^2)+O(\alpha z)
\eq
when $|\alpha z|\leq c$ and $\alpha$ and $c$ are sufficiently small. We have from (\ref{3.10})
\beq \label{3.15}
\frac{1-\psi({-}\alpha z)}{-\alpha^2z/\rho}=(1+\tfrac12\,\sigma_{\alpha}^2z)\Bigl(1+\frac{O(\alpha^2z)+O(\alpha z^2)}{1+\frac12\,\sigma_{\alpha}^2z}\Bigr).
\eq
Now when $z=\mu_0+i\eta\in C_0$ with $\eta\in\dR$, we have by (\ref{3.3})
\begin{eqnarray} \label{3.16}
|1+\tfrac12\,\sigma_{\alpha}^2z|^2 & = & |\tfrac12+\tfrac12\,i\sigma_{\alpha}^2\eta|^2=\tfrac14+\tfrac14\,\sigma_{\alpha}^4\eta^2 \nonumber \\[3mm]
& = & \tfrac14\,\sigma_{\alpha}^4(\mu_0^2+\eta^2)=\tfrac14\,\sigma_{\alpha}^4\,|z|^2.
\end{eqnarray}
Hence
\beq \label{3.17}
\frac{O(\alpha^2z)+O(\alpha z^2)}{1+\frac12\,\sigma_{\alpha}^2z}=O(\alpha^2)+O(\alpha z),
\eq
and this has modulus $\leq\;1/2$ when $\alpha$ is sufficiently small and $|\alpha z|\leq c$ with $c$ determined by the implicit constants in the $O$'s in (\ref{3.10}). This gives (\ref{3.14}).

To finish the proof of (\ref{3.2}), we split the integral in (\ref{3.13}) into the parts $|z|\leq c/\alpha$ and $|z|\geq c/\alpha$. We have by (\ref{3.14})
\begin{align} \label{3.18}
&{\rm log}(\dE\,[e^{-\alpha sW}])  =  \frac{-1}{2\pi i}\,\il_{\scriptsize{\ba{l} z\in C_0, \\ |z|\leq c/\alpha \ea}}\,\frac{s}{z(s-z)}\,{\rm log}(1+\tfrac12\,\sigma_{\alpha}^2z)\,dz \nonumber \\[3.5mm]
&  +~\il_{\scriptsize{\ba{l} z\in C_0, \\ |z|\geq c/\alpha \ea}}\,\frac{s}{z(s-z)}\,(O(\alpha^2)+O(\alpha z))\,dz  - \frac{1}{2\pi i}\,\il_{\scriptsize{\ba{l} z\in C_0, \\ |z|\leq c/\alpha \ea}}\,\frac{s}{z(s-z)}\,{\rm log}\Bigl(\frac{1-\psi({-}\alpha z)}{-\alpha^2z/\rho}\Bigr)\,dz.
\end{align}
For the first integral on the second line of (\ref{3.18}), we have
\begin{eqnarray} \label{3.19}
& \mbox{} & \Bigl|\il_{\scriptsize{\ba{l} z\in C_0, \\ |z|\leq c/\alpha \ea}}\,\frac{s}{z(s-z)}\,(O(\alpha^2)+O(\alpha z))\,dz\Bigr| \nonumber \\[3.5mm]
& & =~O(\alpha^2)+)\Bigl(\il_{\scriptsize{\ba{l} z\in C_0, \\ |z|\leq c/\alpha \ea}}\,\frac{\alpha}{|z|}\,|dz|\Bigr)=O(\alpha^2)+O(\alpha\,{\rm log}(1/\alpha)),
\end{eqnarray}
uniformly in any bounded set of $s$ with ${\rm Re}(s)\geq\frac12\,\mu_0$. For the second integral on the second line of (\ref{3.18}), we use (\ref{3.7}), so that
\beq \label{3.20}
\frac{1}{4\sigma_{\alpha}^2|z|}\leq\Bigl|\frac{1-\psi({-}\alpha z)}{-\alpha^2z/\rho}\Bigr|\leq\frac{2\rho}{\alpha^2|z|}\,,~~~~z\in C_0.
\eq
Therefore, for $z\in C_0$ and $|z|\geq c/\alpha$,
\beq \label{3.21}
-{\rm log}\,|z|+{\rm log}\Bigl(\frac{1}{4\sigma_{\alpha}^2}\Bigr)\leq{\rm log}\Bigl|\frac{1-\psi({-}\alpha z)}{-\alpha^2z/\rho}\Bigr|\leq{\rm log}\Bigl(\frac{2\rho}{c\alpha}\Bigr).
\eq
Using that $(1-\psi({-}\alpha z))/({-}\alpha^2z/\rho)$ lies in the cut plane $\dC\backslash({-}\infty,0]$, so that its argument is between $-\pi$ and $\pi$, we get
\beq \label{3.22}
\Bigl|{\rm log}\Bigl(\frac{1-\psi({-}\alpha z)}{-\alpha^2z/\rho}\Bigr)\Bigr|=O({\rm log}\,|z|)+O({\rm log}(1/\alpha))+O(1).
\eq
From (\ref{3.22}) it follows that the second integral on the second line of (\ref{3.18}) is $O(\alpha\,{\rm log}(1/\alpha))$, uniformly in any bounded set of $s$ with ${\rm Re}(s)\geq\frac12\,\mu_0$.

We conclude that
\begin{eqnarray} \label{3.23}
{\rm log}(\dE\,[e^{-\alpha sW}])  =  \frac{-1}{2\pi i}\,\il_{\scriptsize{\ba{l} z\in C_0, \\ |z|\leq c/\alpha \ea}}\,\frac{s}{z(s-z)}\,{\rm log}(1+\tfrac12\,\sigma_{\alpha}^2z)\,dz  +~O(\alpha\,{\rm log}(1/\alpha)),
\end{eqnarray}
uniformly in any bounded set of $s$ with ${\rm Re}(s)\geq\frac12\,\mu_0$. We extend the integration range in the integral in (\ref{3.23}) to all $z\in C_0$, at the expense of an error $O(\alpha\,{\rm log}(1/\alpha))$ uniformly in any bounded set of $s$ with ${\rm Re}(s)\geq\frac12\,\mu_0$. Finally,
\beq \label{3.24}
\frac{-1}{2\pi i}\,\il_{C_0}\,\frac{s}{z(s-z)}\,{\rm log}(1+\tfrac12\,\sigma_{\alpha}^2z)\,dz={-}{\rm log}(1+\tfrac12\,\sigma_{\alpha}^2s),~~~~~~{\rm Re}(s)\geq\tfrac12\,\mu_0,
\eq
by Cauchy's theorem, and we get (\ref{3.2}).

\section{Kingman-type result for nearly deterministic queues} \label{sec4}

We consider
\beq \label{4.1}
W_n\stackrel{d}{=}\Bigl(W_n+V_n-\frac{1}{\rho_n}\,U_n\Bigr)^+,
\eq
with
\beq \label{4.2}
V_n=\frac1n\,\sum_{k=1}^n\,V_{n,k},~~~~~~U_n=\frac1n\,\sum_{k=1}^n\,U_{n,k},
\eq
where $0<\rho_n<1$, $V_{n,k}$ are i.i.d.\ copies of $V$ and $U_{n,k}$ are i.i.d.\ copies of $U$, with $V$ and $U$ as in Section~\ref{sec2}, and $V_{n,k}$ and $U_{n,k}$ are independent. We shall show that, when $(1-\rho_n)\asymp 1/n$ (so that there are fixed $\beta_1$, $\beta_2$ with $0<\beta_1\leq\beta_2<\infty$ such that $(1-\rho_n)\,n\in[\beta_1,\beta_2]$ for $n=1,2,...$),
\beq \label{4.3}
{\rm log}(\dE\,[e^{-tW_n}])={-}{\rm log}(1+\gamma_nt)+O\Bigl(\frac{{\rm log}\,n}{\sqrt{n}}\Bigr),
\eq
uniformly in any bounded set of $t$ with ${\rm Re}(t)\geq\frac12\,x_0$. Here
\beq \label{4.4}
\gamma_n=\frac{\sigma_V^2+\rho^{-2}\sigma_U^2}{2n(1-\rho_n)}\,\rho_n,~~~~~~x_0=\frac{-1}{2\gamma_n},
\eq
with $\gamma_n$ bounded away from 0 and $\infty$.

We use again Pollaczek's result, so that
\beq \label{4.5}
{\rm log}(\dE\,[e^{-tW_n}])=\frac{-1}{2\pi i}\,\il_C\,\frac{t}{z(t-z)}\,{\rm log}(1-\varp({-}z))\,dz,
\eq
where for $-n\delta\leq{\rm Re}(z)\leq n\delta$ (with $\delta$ as in the first paragraph of Section~\ref{sec2})
\beq
\label{4.6}
\varp({-}z)=\dE\,[e^{-z(V_n-\frac{1}{\rho_n}\,U_n)}]=(\dE\,[e^{-\frac{z}{n}\,(V-\frac{1}{\rho_n}\,U)}])^n=\psi^n({-}z/n),
\eq
with $\psi$ as in (\ref{2.1}) with $\rho=\rho_n$, and $C$ is a line parallel to, and to the left of, the imaginary axis, and to the right of the singularities of ${\rm log}(1-\varp({-}z))$ in the open left half-plane, and lies to the right of $C$.

Assume that $1-\rho_n\asymp 1/n$, and delete $n$ from $\rho_n$ and $\gamma_n$ below for conciseness. We use (\ref{2.3}) with $\zeta=z/n$ and $z=O(\sqrt{n})$, so that
\beq \label{4.7}
\psi({-}z/n)=1+\Bigl(\frac{1}{\rho}-1\Bigr)\,\frac{z}{n}+\Bigl(\frac{1}{\rho}-1\Bigr)\,\gamma\,\frac{z^2}{n}+O\Bigl(\frac{z^2}{n^4}\Bigr)+O\Bigl(\frac{z^3}{n^3}\Bigr).
\eq
With $\varp({-}z)$ as in \eqref{4.6} and using the expansion $(1+X)^n=1+nX+\frac12\,n^2X^2+O(n^3X^3)$, valid for $X=O(\frac1n)$, we get
\beq \label{4.8}
\varp({-}z)=1+\Bigl(\frac{1}{\rho}-1\Bigr)\,z+\Bigl(\frac{1}{\rho}-1\Bigr)\,\gamma z^2+O\Bigl(\frac{|z|^2+|z|^3+|z|^4}{n^2}\Bigr).
\eq
In the $O$ in (\ref{4.8}) terms like $z^6/n^4$, which is $O(z^2/n^2)$ when $z=O(\sqrt{n})$, have been collected. The leading part of the right-hand side of (\ref{4.8}), considered for $z\leq0$, equals 1 for $z=0$ and $z={-}1/\gamma$, and is minimal for $z={-}1/2\gamma=x_0$, with minimum value $1-(1-\rho)/4\rho\gamma$. Therefore, for large $n$,
\beq \label{4.9}
|\varp({-}x_0-iy)|\leq\varp({-}x_0)\leq1-(1-\rho)/8\gamma,~~~~~~y\in\dR.
\eq
Hence, we can use $C=\{x_0+iy\:|\:y\in\dR\}$ in (\ref{4.5}), so that $1-\varp({-}z)$ has positive real part when $z\in C$ and $n$ is large, with principal-value logarithm for the log in the integral.

We have for $z\in C$, $z=O(\sqrt{n})$ from (\ref{4.8}) and $(1-\rho_n)\asymp1/n$,
\beq \label{4.10}
\frac{1-\varp({-}z)}{-(\frac{1}{\rho}-1)\,z}=1+\gamma z+O\Bigl(\frac{|z|+|z|^2+|z|^3}{n}\Bigr).
\eq
We are now in a quite similar position as in Section~\ref{sec3} from (\ref{3.10}) onwards. In particular, using $|1+\gamma z|=\gamma|z|$ for $z\in C$, we have
\begin{eqnarray} \label{4.11}
{\rm log}\Bigl(\frac{1-\varp({-}z)}{-(\frac{1}{\rho}-1)\,z}\Bigr) & = & {\rm log}(1+\gamma z)+O\Bigl(\frac{|z|+|z|^2+|z|^3}{n(1+\gamma z)}\Bigr) \nonumber \\[3.5mm]
& = & {\rm log}(1+\gamma z)+O\Bigl(\frac{1+|z|+|z|^2}{n}\Bigr),
\end{eqnarray}
when $|z|\leq c\sqrt{n}$ and $c>0$ is small enough to ensure that the $O$-terms in (\ref{4.11}) do not exceed $1/2$. Furthermore,
\begin{eqnarray} \label{4.12}
& \mbox{} & \il_{\scriptsize{\ba{l} z\in C, \\ |z|\leq c\sqrt{n}\ea}}\,\Bigl|\frac{t}{z(t-z)}\Bigr|\,\frac{1+|z|+|z|^2}{n}\,|dz| \nonumber \\[3.5mm]
& & =O\Bigl(\il_{\scriptsize{\ba{l} z\in C, \\ |z|\leq c\sqrt{n}\ea}}\,\frac{1+|z|+|z|^2}{n\,|z|^2}\,|dz|\Bigr)=O\Bigl(\frac{1}{\sqrt{n}}\Bigr),
\end{eqnarray}
uniformly in any bounded set of $t$ with ${\rm Re}(t)\geq\frac12\,x_0$. Proceeding then as in Section~\ref{sec3} from (\ref{3.18}) onwards, with (\ref{4.12}) as substitute for (\ref{3.19}) and $\sqrt{n}$ instead of $1/\alpha$, we get (\ref{4.3}).

\section{Gaussian regime for nearly deterministic queues} \label{sec5}

We consider the same setting as in Section~\ref{sec4}, but now we assume that $(1-\rho_n)\asymp 1/\sqrt{n}$, so that there are fixed $\beta_1$, $\beta_2$ with $0<\beta_1\leq\beta_2<\infty$ such that $(1-\rho_n)\,\sqrt{n}\in [\beta_1,\beta_2]$ for $n=1,2,...\,$. The precise formulation of our results requires quantities derived from $h(\zeta)={\rm log}\,(\psi({-}\zeta)$ at the saddle point $\zeta_{sp}$ of $h$ on the negative real axis, and an additional condition discussed below.

We shall show the following (Theorem~\ref{thm3}). Let
\beq \label{5.1}
\sigma_n=(h''(\zeta_{sp}))^{1/2},~~~~~~\beta_n={-}\zeta_{sp}\,\sigma_n\,\sqrt{n}.
\eq
Then we have $\sigma_n\asymp 1$, $\beta_n\asymp 1$ as $n\pr\infty$, and
\beq \label{5.2}
{\rm log}\Bigl(\dE\,\Bigl[\exp\Bigl({-}s\,\frac{\sqrt{n}}{\sigma_n}\,W_n\Bigl)\Bigr]\Bigr)={\rm log}(\dE\,[e^{-sM_{\beta_n}}])+O\Bigl(\frac{1}{\sqrt{n}}\Bigr),~~~~~~n\pr\infty,
\eq
uniformly in any bounded set of $s$ with ${\rm Re}(s)\geq{-}\frac12\,\beta_n$. We shall also show the following refinement of Theorem~\ref{thm3}.

\begin{thm} \label{thm4}
Let $\sigma_n$ and $\beta_n$ as in (\ref{5.1}), and let
\beq \label{5.3}
d_2={-}\,\frac{h'''(\zeta_{sp})}{6h''(\zeta_{sp})}.
\eq
Then $d_2=O(1)$, $n=1,2,...\,$, and, uniformly in any bounded set of $s$ with ${\rm Re}(s)\geq{-}\,\frac12\,\beta_n$,
\beq \label{5.4}
{\rm log}\Bigl(\dE\,\Bigl[\exp\Bigl({-}s\,\frac{\sqrt{n}}{\sigma_n}\,W_n\Bigr)\Bigr]\Bigr)={\rm log}(\dE\,[e^{-R_nM_{B_n}}])+O\Bigl(\frac1n\Bigr),~~~~~~n\pr\infty,
\eq
where
\beq \label{5.5}
B_n=\frac{\beta_n}{1+\beta_n\varp_n},~~~~~~R_n=R_n(s)=\frac{s}{(1+\beta_n\varp_n)(1+(s+\beta_n)\varp_n)},
\eq
with $\varp_n=d_2/\sigma_n\sqrt{n}=O(1/\sqrt{n})$.
\end{thm}

%From (\ref{5.2}) the moments of $\frac{\sqrt{n}}{\sigma_n}\,W$ are approximated with error $O(1/\sqrt{n})$ by the moments $m_k(\beta_n)$ of $M_{\beta_n}$ as in (\ref{1.11}). Also from the refinement Theorem~\ref{thm4}, we have for $k=1,2,...$
%\beq \label{5.6}
%\dE\,\Bigl[\Bigl(\frac{\sqrt{n}}{\sigma_n}\,W_n\Bigr)^k\Bigr]=({-}1)^k\,\Bigl(\frac{d}{ds}\Bigr)^k\,(\dE\,[e^{-R_n(s)M_{B_n}}])\Bigl|_{s=0}+O\Bigl(\frac1n\Bigr).
%\eq
%This yields for $k=1,2,...$
%\beq \label{5.7}
%\dE\,\Bigl[\Bigl(\frac{\sqrt{n}}{\sigma_n}\,W_n\Bigr)^k\Bigr]=\frac{m_k(B_n)}{(1+\beta_n\varp_n)^{2k}}+\frac{k(k-1)m_{k-1}(B_n)\varp_n}{(1+\beta_n\varp_n)^{2k-1}}+O\Bigl(\frac1n\Bigr),
%\eq
%with $\varp_n$ as in Theorem~\ref{thm4} and $m_k(B_n)=\dE\,[M_{B_n}^k]$.

From (\ref{5.2}) the moments of $\frac{\sqrt{n}}{\sigma_n}\,W_n$ are approximated with error $O(1/\sqrt{n})$ by the moments $m_k(\beta_n)$ of $M_{\beta_n}$ as in (\ref{1.11}).

As a consequence of Theorem~\ref{thm4}, we have the following result.
\begin{thm}\label{thm5}
For $k=1,2,...$
\begin{align}
\dE\,\Bigl[\Bigl(\frac{\sqrt{n}}{\sigma_n}\,W_n\Bigr)^k\Bigr]\,&\,=({-}1)^k\,\Bigl(\frac{d}{ds}\Bigr)^k\,(\dE\,[e^{-R_n(s)M_{B_n}}])\Bigl|_{s=0}+O\Bigl(\frac1n\Bigr)
\nonumber\\
&\,=\frac{m_k(B_n)}{(1+\beta_n\varp_n)^{2k}}+\frac{k(k-1)m_{k-1}(B_n)\varp_n}{(1+\beta_n\varp_n)^{2k-1}}+O\Bigl(\frac1n\Bigr),
\label{5.6}
\end{align}
with $\varp_n$ as in Theorem~\ref{thm4} and $m_k(B_n)=\dE\,[M_{B_n}^k]$.
\end{thm}

%This yields for $k=1,2,...$
%\beq \label{5.7}
%\dE\,\Bigl[\Bigl(\frac{\sqrt{n}}{\sigma_n}\,W_n\Bigr)^k\Bigr]=\frac{m_k(B_n)}{(1+\beta_n\varp_n)^{2k}}+\frac{k(k-1)m_{k-1}(B_n)\varp_n}{(1+\beta_n\varp_n)^{2k-1}}+O\Bigl(\frac1n\Bigr),
%\eq
%with $\varp_n$ as in Theorem~\ref{thm4} and $m_k(B_n)=\dE\,[M_{B_n}^k]$.

We shall prove Theorem~\ref{thm3} below, and we present the proofs of Theorem~\ref{thm4} and Theorem~\ref{thm5} in Appendices~A and B. For all these proofs, we use again Pollaczek's formula (\ref{4.5}) in which we substitute $\zeta=z/n$. Thus, we have
\beq \label{5.8}
{\rm log}(\dE\,[e^{-tW_n}])=\frac{-1}{2\pi i}\,\il_{C_n}\,\frac{t/n}{\zeta(t/n-\zeta)}\,{\rm log}(1-\psi^n({-}\zeta))\,d\zeta,
\eq
where $C_n=\frac1n\,C$ is a line parallel to, and to the left of, the imaginary axis, and to the right of the singularities of ${\rm log}(1-\psi^n({-}\zeta))$, and $t/n$ lies to the right of $C_n$.
We now need the following assumption, because it allows us to integrate in (\ref{5.8}) over $\zeta$ with $|{\rm Im}(\zeta)|\leq\delta$, at the expense of an error of exponential decay as $n\pr\infty$.

\begin{assumption}\label{ass:delta}
There is a $\delta>0$, $\eps>0$ such that for all  $\rho\in[\frac12,1)$ and $x\in[{-}\delta,0)$ and $y\in\dR$, $|y|>\delta$, we have $|\psi({-}x+iy)|<\psi({-}x)-\eps$.
\end{assumption}

With reference to Section~\ref{sec2}, we can assume that $\delta>0$ is such that $h(\zeta)={\rm log}\,\psi({-}\zeta)$ is analytic in the rectangle $-\delta\leq{\rm Re}(\zeta)\leq0$, $|{\rm Im}(\zeta)|\leq\delta$, and, by taking $n$ sufficiently large with $(1-\rho)=(1-\rho_n)\asymp 1/\sqrt{n}$ in (\ref{2.7}), that the saddle point $\zeta_{sp}$ of $h$ lies in this rectangle. We then have, with exponentially small error as $n\pr\infty$,
\beq \label{5.9}
{\rm log}(\dE\,[e^{-t W_n}])=\frac{-1}{2\pi i}\,\il_{\zeta_{sp}-i\delta}^{\zeta_{sp}+i\delta}\,\frac{t/n}{\zeta(t/n-\zeta)}\,{\rm log}(1-e^{nh(\zeta)})\,d\zeta
\eq
when ${\rm Re}(t/n)\geq\frac12\,\zeta_{sp}$.

We have from (\ref{2.6}--\ref{2.8}) and $(1-\rho_n)\asymp 1/\sqrt{n}$ that
\beq \label{5.10}
\zeta_{sp}\asymp\frac{1}{\sqrt{n}},~~~~~~h(\zeta_{sp})\asymp\frac1n,~~~~~~\sigma_n^2=h''(\zeta_{sp})\asymp1,
\eq
and this shows that $\sigma_n\asymp1$, $\beta_n\asymp1$ and also that $d_2=O(1)$, see (\ref{5.1}) and (\ref{5.2}).

For both Theorem~\ref{thm3} and \ref{thm4}, we shall work from the integral in (\ref{5.9}) towards the integral representation of the mgf $\dE\,[\exp({-}s\,M_{\beta})]$ of the maximum $M_{\beta}$ of the Gaussian random walk with drift $-\beta$. For the latter we have (from Pollaczek's formula, applied to $W=(W+V-U)^+$ with $V$ and $U$ having pdf's $\frac{1}{\sqrt{2\pi}}\,\exp({-}\frac12\,(x-A)^2)\,\chi_{[0,\infty)}(x)$ and $\delta_{\beta}(x-A)$, respectively, and letting $A\pr\infty$)
\beq \label{5.11}
{\rm log}(\dE\,[e^{-s M_{\beta}}])=\frac{-1}{2\pi i}\,\il_C\,\frac{s}{z(s-z)}\,{\rm log}(1-e^{\beta z+\frac12 z^2})\,dz,
\eq
where $C$ is a line parallel to, and to the left of, the imaginary axis, and $s\in C$ is to the right of $C$. Substituting $z={-}\beta+iy$, $-\infty<y<\infty$, we get for $s\in\dC$, ${\rm Re}(s)>{-}\beta$
\beq \label{5.12}
{\rm log}(\dE\,[e^{-sM_{\beta}}])=\frac{1}{2\pi}\,\il_{-\infty}^{\infty}\,\frac{s}{(\beta-iy)(s+\beta-iy)}\,{\rm log}(1-e^{-\frac12\beta^2-\frac12 y^2})\,dy.
\eq

We make in the integral in (\ref{5.9}) a substitution that brings $\exp(nh(\zeta))$ in Gaussian form. We thus let $\zeta=\zeta(v)$ be the solution of the equation
\beq \label{5.13}
h(\zeta(v))=h(\zeta_{sp})-\tfrac12\,v^2\,h''(\zeta_{sp})
\eq
that satisfies $\zeta(v)=\zeta_{sp}+iv+O(v^2)$ as $v\pr0$. By Lagrange's theorem, there is an $r>0$, independent of $n$, such that
\beq \label{5.14}
\zeta(v)=\zeta_{sp}+iv+\sum_{l=2}^{\infty}\,d_l(iv)^l,~~~~~~|v|\leq r,
\eq
with real $d_l$, and $d_2$ given by (\ref{5.3}), see \cite{ref9}, end of Section~3 for more details about such a substitution. With the substitution $\zeta=\zeta(v)$, $-r\leq v\leq r$, in (\ref{5.9}), we get with exponentially small error
\beq \label{5.15}
{\rm log}(\dE\,[e^{-tW_n}])=\frac{-1}{2\pi i}\,\il_{-r}^r\,\frac{\zeta'(v)\,t/n}{\zeta(v)(t/n-\zeta(v))}\,{\rm log}(1-e^{nh(\zeta_{sp})-\frac12 v^2h''(\zeta_{sp})})\,dv
\eq
when $n\pr\infty$ and ${\rm Re}(t)\geq\frac12\,n\zeta_{sp}$. We write for conciseness $\sigma=\sigma_n$ in the sequel, and we take $t=s\,\sqrt{n}/\sigma$ in (\ref{5.16}) and substitute $y=v\sigma\sqrt{n}$, see (\ref{5.1}). Then, for $s$ in a bounded set contained in ${\rm Re}(s)\geq{-}\frac12\,\beta_n$ (so that ${\rm Re}(t)\geq{-}\frac12\,\beta_n\,\sqrt{n}/\sigma=\frac12\,n\,\zeta_{sp}$, see (\ref{5.1})), we have with exponentially small error
\begin{eqnarray} \label{5.16}
& \mbox{} & {\rm log}\Bigl(\dE\,\Bigl[\exp\Bigl({-}\,\frac{s\sqrt{n}}{\sigma}\,W_n\Bigr)\Bigr]\Bigr) \nonumber \\[3.5mm]
& & =~\frac{1}{2\pi}\,\il_{-R}^R\,\frac{is\zeta'(y/\sigma\sqrt{n})\,{\rm log}(1-e^{nh(\zeta_{sp})-\frac12 y^2})} {\sigma\sqrt{n}\,\zeta(y/\sigma\sqrt{n})(s-\sigma\sqrt{n}\,\zeta(y/\sigma\sqrt{n}))}\,dy,
\end{eqnarray}
where $R=r\sigma\sqrt{n}$.

The remainder of the proofs of (\ref{5.2}) and Theorem~\ref{thm4} consists now of approximating $nh(\zeta_{sp})$ by $-\frac12\,\beta_n^2$ and $-\frac12\,B_n^2$, respectively, and approximating the front factor
\beq \label{5.17}
FF=\frac{is\zeta'(y/\sigma\sqrt{n})}{\sigma\sqrt{n}\,\zeta(y/\sigma\sqrt{n})(s-\sigma\sqrt{n}\,\zeta(y/\sigma\sqrt{n}))},
\eq
using a linear and quadratic approximation, respectively, from the power series in (\ref{5.14}).

For (\ref{5.2}) we thus use in (\ref{5.16})
\beq \label{5.18}
\zeta'(y/\sigma\sqrt{n})=i+O(y/\sqrt{n})\,,~~~~~\sigma\sqrt{n}\,\zeta(y/\sigma\sqrt{n})={-}\beta_n+iy+O(y^2/\sqrt{n}),
\eq
and obtain, uniformly in any bounded set of $s$ such that ${\rm Re}(s)\geq{-}\frac12\,\beta_n$,
\begin{eqnarray} \label{5.19}
& \mbox{} & {\rm log}\Bigl(\dE\,\Bigl[\exp\Bigl({-}\,\frac{s\sqrt{n}}{\sigma_n}\,W_n\Bigr)\Bigr]\Bigr) \nonumber \\[3.5mm]
& & =~\frac{1}{2\pi}\,\il_{-R}^R\,\frac{s\,{\rm log}(1-e^{-\frac12\beta_n^2-\frac12 y^2})}{(\beta_n-iy)(s+\beta_n-iy)}\,dy\Bigl(1+O\Bigl(\frac{1}{\sqrt{n}}\Bigr)\Bigr),
\end{eqnarray}
where we have restored $\sigma=\sigma_n$ on the left-hand side of (\ref{5.19}). Recalling that $R=r\,\sigma_n\sqrt{n}$, with $\sigma_n\asymp1$ and $r>0$ independent of $n$, we see that the integral in (\ref{5.19}) equals the integral in (\ref{5.12}), with $\beta_n$ instead of $\beta$, within exponentially small error when ${\rm Re}(s)\geq{-}\frac12\,\beta_n$. This completes the proof of Theorem~\ref{thm3}.

The proof of Theorem~\ref{thm4} is similar, but the details require somewhat more elaboration, and are therefore given in Appendix~\ref{appA}. We show Theorem~\ref{thm5} in Appendix~\ref{appB} and also that the $O(1/\sqrt{n})$ at the right-hand side of \eqref{1.10} can be replaced by $O(1/n)$ when the third cumulants of $V$ and $U$ are equal.

\section{Computational issues} \label{sec6}

In this section we present computational schemes for
exact and approximate values of the (properly scaled) moments of the steady-state waiting time distribution $W$ via Pollaczek's formula \eqref{1.3},
\beq\label{75}
\log\big(\dE[e^{-sW}]\big)=\frac{-1}{2\pi i}\int_C \frac{s}{z(s-z)}\log(1-\psi(-z))\,dz,
\eq
with $\psi(-z)=\dE[\exp(-z(V-\frac1\rho U))]$, as earlier. In Section~\ref{sec7} this is made specific by making various choices for $V$ and $U$.

\subsection{Moments and cumulants}\label{subsec61}

Assume that $X$ is a random variable with finite moments of all order. We compute the moments
\beq\label{76}
m_k(X)=\dE[X^k]=\left(\frac{d}{ds}\right)^k\left(\dE[e^{sX}]\right)_{s=0}\, , k=0,1,\dots,
\eq
of $X$ from the cumulants
\beq\label{77}
c_l(X)=\left(\frac{d}{ds}\right)^l\log\left(\dE[e^{sX}]\right)_{s=0}\, , l=0,1,\dots,
\eq
of $X$ recursively according to
\beq\label{78}
m_0(X)=1 \,;\, m_k(X)=\sum_{l=1}^k\binom{k-1}{l-1}c_l(X)m_{k-l}(X), \, k=1,2,\dots \, .
\eq

\subsection{Moments in the classical heavy-traffic Kingman case}\label{subsec62}

With $\alpha=1-\rho\,\downarrow\,0$, we compute the moments $m_k(\alpha W)=\dE[(\alpha W)^k]$ of $\alpha W$ from the cumumants of $\alpha W$ according to \eqref{78} in Subsection~\ref{subsec61}. The latter are obtained from the appropriate version \eqref{3.8} of Pollaczek's formula, so that
\begin{align}
c_l(\alpha W)\,&=\,(-1)^l\left(\frac{d}{ds}\right)^l\left(\frac{-1}{2\pi i}\int_{C_0}\frac{z}{z(s-z)}\log(1-\psi(-\alpha z))dz\right)_{s=0}\nonumber\\
&=\,\frac{(-1)^l l!}{\pi}\int_0^\infty \mbox{Re}\left[\frac{\log(1-\psi(-\alpha z))}{z^{l+1}}\right]dy,
\label{79}
\end{align}
where $C_0=\{z=\mu_0+iy\,|\,-\infty<y<\infty\}$ with $\mu_0=-1/\sigma_\alpha^2=-1/[(\sigma_V^2+\rho^{-2}\sigma^2_U)\rho]$, compare \cite[Equation (15)]{ref1}.

From the result of Section~\ref{sec3}, see \eqref{1.7} in Theorem \ref{thm1}, we have for $k=1,2,\dots$
\begin{equation}\label{eqn80}
\dE[(\alpha W)^k]=k! (\frac12 \sigma_\alpha^2)^k+O(\alpha\log(1/\alpha)), \quad \alpha\downarrow0,
\end{equation}
where $\sigma_\alpha^2=(\sigma_V^2+\rho^{-2}\sigma^{2}_U)\rho$.

\subsection{Moments in the nearly deterministic heavy-traffic Kingman case}\label{subsec63}

With $1-\rho_n \asymp 1/n$ and $W_n$ as in Section~\ref{sec4}, we compute the moments $m_k(W_n)=\dE[W_n^k]$ of $W_n$ from the cumulants $c_l(W_n)$ of $W_n$ according to \eqref{78} in Subsection~\ref{subsec61}. The latter are obtained from the appropriate version \eqref{4.5} of Pollaczek's formula, so that
\begin{align}
c_l(W_n)\,&=\,(-1)^l\left(\frac{d}{ds}\right)^l\left(\frac{-1}{2\pi i}\int_{C}\frac{z}{z(s-z)}\log(1-\psi^n(-z/n))dz\right)_{s=0}\nonumber\\
&=\,\frac{(-1)^l l!}{\pi}\int_0^\infty \mbox{Re}\left[\frac{\log(1-\psi^n(-z/n))}{z^{l+1}}\right]dy,
\label{81}
\end{align}
where $C=\{z=\frac{-1}{2\gamma_n}+iy\,|\,-\infty<y<\infty\}$ and
\begin{equation}
\label{82}
\gamma_n=(\sigma_V^2+\rho_n^{-2}\sigma_U^2)\rho_n/(2n(1-\rho_n)).
\end{equation}

From the result of Section~\ref{sec4}, see \eqref{1.9} in Theorem~\ref{thm2}, we have for $k=1,2,\dots$
\begin{equation}\label{eqn83}
\dE[W_n^{k}]=k!\gamma_n^k+O\left(\frac{\log n}{\sqrt{n}}\right), \quad n\to\infty.
\end{equation}

\subsection{Moments in the nearly deterministic heavy-traffic Gaussian case}\label{subsec64}

We consider the dedicated saddle point method with $1-\rho_n \asymp 1/\sqrt{n}$ and $W_n$ as given in Section~\ref{sec5}. Thus $\zeta_{sp}$ is the unique saddle point $\zeta$ on the negative real axis of $h(\zeta)=\log \psi(-\zeta)$, characterized by $h'(\zeta)=\psi'(-\zeta)=0$, see Section~\ref{sec2}. In several of the specific cases as given in Section~\ref{sec7}, $\zeta_{sp}$ can be found in closed form; in general a Newton iteration can be used, using the leading term at the right-hand side of \eqref{2.7} as initial value. We then let
\begin{equation}\label{eqn84}
\sigma_n=(h''(\zeta_{sp}))^{1/2}, \ \beta_n=-\zeta_{sp}\sigma_n\sqrt{n}.
\end{equation}

The moments $m_k\left(\frac{\sqrt{n}}{\sigma_n}W_n\right)=\dE\left[\left(\frac{\sqrt{n}}{\sigma_n}W_n\right)^k\right]$ of $\frac{\sqrt{n}}{\sigma_n}W_n$ can be computed from the cumulants $c_l\left(\frac{\sqrt{n}}{\sigma_n}W_n\right)$ according to \eqref{78} in Subsection~\ref{subsec61}. The latter are obtained from the appropriate version \eqref{5.8} of Pollaczek's formula, with $t=s\sqrt{n}/\sigma_n$, so that
\begin{align}
c_l\left(\frac{\sqrt{n}}{\sigma_n} W_n\right)\,&=\,(-1)^l\left(\frac{d}{ds}\right)^l\left(\frac{-1}{2\pi i}\int_{C_n}\frac{(s/\sigma_n\sqrt{n})\log(1-\psi^n(-\zeta))}{\zeta((s/\sigma_n\sqrt{n})-\zeta)}d\zeta\right)_{s=0}\nonumber\\
&=\,\frac{(-1)^l l!}{\pi}\left(\frac{1}{\sigma_n\sqrt{n}}\right)^l\int_0^\infty \mbox{Re}\left[\frac{\log(1-\psi^n(-\zeta))}{\zeta^{l+1}}\right]dy,
\label{eqn85}
\end{align}
where $C_n=\{\zeta=\zeta_{sp}+iy\,|\,-\infty<y<\infty\}$.

From the result in Section~\ref{sec5}, see \eqref{1.11} in Theorem~\ref{thm3}, we have for $k=1,2,\dots$
\begin{align}\label{eqn86}
\dE\left[\left(\frac{\sqrt{n}}{\sigma_n}W_n\right)^k\right]=m_k(\beta_n)+O\left(\frac{1}{\sqrt{n}}\right), \quad n\to\infty,
\end{align}
where $m_k(\beta)$ is the $k^\text{th}$ moment of the maximum $M_\beta$ of the Gaussian random walk with drift $-\beta$. These $m_k(\beta)$ can be computed from the cumulants $c_l(M_\beta)$ of $M_\beta$ using \eqref{78} in Subsection~\ref{subsec61}. The latter can be computed by numerical integration, using \eqref{5.12}, so that
\begin{align}
c_l\left(M_\beta\right)\,&=\,\frac{(-1)^l l!}{\pi}\int_0^\infty \mbox{Re}\left[\frac{\log(1-e^{\beta z+\frac12 z^2})}{z^{l+1}}\right]dy,
\label{eqn87}
\end{align}
where $z=-\beta+iy, -\infty<y<\infty$. Alternatively, we have from Theorem~\ref{thm1} in \cite{ref8} for $0<\beta<2\sqrt{\pi}$ and $l=1,2,\dots$
\begin{align}\label{eqn88}
c_l\left(M_\beta\right)\,&=\,\frac{(l-1)!}{(2\beta)^l}+\frac{1}{\sqrt{2\pi}}\sum_{j=0}^l\binom{l}{j}\Gamma\left(\frac{l-j+1}{2}\right) \zeta\left(-\frac12 l-\frac12 j+1\right)2^{\frac{l-j-1}{2}}(-\beta)^j\nonumber\\
&+\,\frac{(-1)^{l+1}l!}{\sqrt{2\pi}}\sum_{r=0}^{\infty}\frac{\zeta(-l-r+\frac12)(-\frac12)^r\beta^{2r+l+1}}{r!(2r+1)\cdot\ldots\cdot(2r+l+1)},
\end{align}
where $\zeta$ denotes here the Riemann zeta function (not to be confused with the function $\zeta(v)$ in (\ref{5.13}--\ref{5.14}) pertaining to the dedicated saddle point method).

From Theorems \ref{thm4} and \ref{thm5} in Section~\ref{sec5}, we can refine the result in \eqref{eqn86} according to
\be\label{eqn89}
\dE\,\Bigl[\Bigl(\frac{\sqrt{n}}{\sigma_n}\,W_n\Bigr)^k\Bigr]=\frac{m_k(B_n)}{(1+\beta_n\varp_n)^{2k}}+\frac{k(k-1)\,m_{k-1}(B_n)\varp_n}{(1+\beta_n\varp_n)^{2k-1}}+O\Bigl(\frac1n\Bigr)~,
\eq
where
\be\label{eqn90}
B_n=\beta_n/(1+\beta_n\varp_n), \qquad \varp_n=-\frac{h'''(\zeta_{sp})}{6h''(\zeta_{sp})\sigma_n\sqrt{n}},
\eq
with $\sigma_n$ and $\beta_n$ given in \eqref{eqn84}.

In several of the specific cases as mentioned in Section~\ref{sec7}, the quantities $h''(\zeta_{sp})$ and $h'''(\zeta_{sp})$ as required in \eqref{eqn84} and \eqref{eqn90} can be computed in closed form. In general, one has
\begin{align}
h''(\zeta_{sp})\,&=\,\frac{f_V''}{f_V}-\left(\frac{f_V'}{f_V}\right)^2+\left(\frac{f_U''}{f_U}-\left(\frac{f_U'}{f_U}\right)^2\right) \frac{1}{\rho^2},\label{eqn91}\\
h'''(\zeta_{sp})\,&=\,-\frac{f_V'''}{f_V}+3\frac{f_V''f_V'}{f_V^2}+\left(\frac{f_U'''}{f_U}-3\frac{f_U''f_U'}{f_U^2}\right) \frac{1}{\rho^3},\label{eqn92}
\end{align}
where $f_V(z)=\dE[\exp(zV)]$, $f_U(z)=\dE[\exp(zU)]$, and where all $f_V^{(l})$ and $f_U^{(l)}$ appearing at the right-hand sides of (\ref{eqn91}--\ref{eqn92}) have to be evaluated at $z=-\zeta_{sp}$ and $z=\zeta_{sp}/\rho$, respectively.

We end this section with a  note on the computational issues encountered when evaluating the contour integrals in expressions such as \eqref{79},\eqref{81},\eqref{eqn85} and \eqref{eqn87}. Although modern computer algebra software can numerically evaluate these integrals, despite one of the limits being infinity, one has to carefully choose a sufficiently high numerical accuracy in order to obtain accurate results. Not all software packages may support this feature, which is the reason why we used Wolfram Mathematica for the numerical computations in this paper. Unsurprisingly, this may lead to long computation times, in particular when the load of the system is close to one. In contrast, the numerical evaluation of the approximations in this paper takes practically zero time and does not suffer from numerical issues.

\subsection{Example 1: $U$ and $V$ are Gamma distributed}

We now demonstrate the results for a specific example. Further numerical examples are provided Appendix~\ref{sec7}. In this first example, we consider the case that both $V$ and $U$ have a Gamma distribution, with means $k_U\vart_U=k_V\vart_V=1$, variances $\sigma^2_V=k_V\vart^2_V=\vart_V$ and $\sigma^2_U=k_U\vart^2_U=\vart_U$, and pdfs
\beq \label{6.1}
\frac{x^{k_V-1}\,e^{-x/\vart_V}}{\Gamma(k_V)\,\vart_V^{k_V}}\,,~~x\geq0~;~~~~~~\frac{x^{k_U-1}\,e^{-x/\vart_U}}{\Gamma(k_U)\,\vart_U^{k_U}}\,,~~x\geq0~,
\eq
respectively, with $0<\vart_V,\vart_U<\infty$. Note that Assumption~\ref{ass:delta} in Section~\ref{sec5} is easily checked.
The third cumulants of $U$ and $V$ are $c_3(U)=2\vart_U^2$ and $c_3(V)=2\vart_V^2$, respectively.
Then
%\beq \label{6.2}
%\dE\,[e^{-\zeta V}]=(1+\vart_V\zeta)^{-k_V}~,~~~~~~{\rm Re}(\zeta)>{-}\vart_V^{-1}~,
%\eq
%\beq \label{6.3}
%\dE\,[e^{\zeta U/\rho}]=(1-\vart_U\zeta/\rho)^{-k_U}~,~~~~~~{\rm Re}(\zeta)<\rho\,\vart_U^{-1}~,
%\eq
%and
\beq \label{6.4}
\psi({-}\zeta)=\dE\,[e^{-\zeta(V-\frac{1}{\rho}U)}]=(1+\vart_V\zeta)^{-k_V}\,(1-\vart_U\zeta/\rho)^{-k_U}
\eq
for $-\vart_V^{-1}<{\rm Re}(\zeta)<\rho\,\vart_U^{-1}$. As a consequence of $U$ and $V$ both being Gamma distributed, all approximations can be obtained in closed-form. In our numerical experiments we take three different parameter sets $(\vart_U, \vart_V)$ = $(5/2, 1/2)$, $(1/2, 5/2)$ and $(3/2, 3/2)$.

\subsubsection{Classical regime} We approximate the moments of $\alpha W$, for $\alpha\downarrow0$, by \eqref{eqn80}, with
\beq \label{6.9}
\sigma_{\alpha}^2=(\sigma_V^2+\rho^{-2}\sigma_U^2)\,\rho=(k_V\vart_V^2+\rho^{-2}k_U\vart_U^2)\,\rho=(\vart_V+\rho^{-2}\vart_U)\,\rho~,
\eq
and $\rho=1-\alpha$, where we take $\alpha=1/10,1/100,1/1000$. The results, for the first five moments of $W$, are presented in Table~\ref{tbl:1}.

\begin{table}[t]
\[
\begin{array}{|l|rr|rr|rr|}
\hline
\mc{7}{|c|}{%\text{Example 1  - Classical HT Kingman}~~ (
\vart_U=5/2, \vart_V=1/2}\\
\hline
& \mc{2}{c|}{\alpha=1/10} & \mc{2}{c|}{\alpha=1/100} & \mc{2}{c|}{\alpha=1/1000} \\
\hline
k & {\rm Exact} & {\rm Asymp} & {\rm Exact} & {\rm Asymp} & {\rm Exact} & {\rm Asymp} \\
\hline
 1 & 1.396 & 1.614 & 1.490 & 1.510 & 1.499 & 1.501 \\
 2 & 4.092 & 5.209 & 4.459 & 4.561 & 4.496 & 4.506 \\
 3 & 17.987 & 25.222 & 20.021 & 20.663 & 20.227 & 20.291 \\
 4 & 105.426 & 162.819 & 119.860 & 124.814 & 121.336 & 121.825 \\
 5 & 772.403 & 1313.861 & 896.946 & 942.427 & 909.816 & 914.295 \\
\hline
\end{array}
\]
\[
\begin{array}{|l|rr|rr|rr|}
\hline
\mc{7}{|c|}{%\text{Example 1  - Classical HT Kingman}~~ (
\vart_U=1/2, \vart_V=5/2}\\
\hline
& \mc{2}{c|}{\alpha=1/10} & \mc{2}{c|}{\alpha=1/100} & \mc{2}{c|}{\alpha=1/1000} \\
\hline
k & {\rm Exact} & {\rm Asymp} & {\rm Exact} & {\rm Asymp} & {\rm Exact} & {\rm Asymp} \\
\hline
 1 & 1.331 & 1.403 & 1.483 & 1.490 & 1.498 & 1.499 \\
 2 & 4.083 & 3.936 & 4.459 & 4.440 & 4.496 & 4.494 \\
 3 & 18.806 & 16.562 & 20.111 & 19.849 & 20.236 & 20.210 \\
 4 & 115.505 & 92.932 & 120.933 & 118.300 & 121.444 & 121.176 \\
 5 & 886.802 & 651.817 & 909.028 & 881.352 & 911.030 & 908.217 \\
\hline
\end{array}
\]
%\caption{Same as Table~\ref{tbl:1-1}, except that now $\vart_U=1/2$, $\vart_V=5/2$.}
%\label{tbl:1-2}
%\end{table}

%\begin{table}[t]
\[
\begin{array}{|l|rr|rr|rr|}
\hline
\mc{7}{|c|}{%\text{Example 1  - Classical HT Kingman}~~
\vart_U=3/2, \vart_V=3/2}\\
\hline
& \mc{2}{c|}{\alpha=1/10} & \mc{2}{c|}{\alpha=1/100} & \mc{2}{c|}{\alpha=1/1000} \\
\hline
k & {\rm Exact} & {\rm Asymp} & {\rm Exact} & {\rm Asymp} & {\rm Exact} & {\rm Asymp} \\
\hline
 1 & 1.367 & 1.508 & 1.487 & 1.500 & 1.499 & 1.500 \\
 2 & 4.100 & 4.550 & 4.460 & 4.500 & 4.496 & 4.500 \\
 3 & 18.452 & 20.589 & 20.071 & 20.253 & 20.232 & 20.250 \\
 4 & 110.711 & 124.223 & 120.427 & 121.525 & 121.393 & 121.500 \\
 5 & 830.332 & 936.845 & 903.203 & 911.480 & 910.446 & 911.252 \\
\hline
\end{array}
\]
\caption{Example 1  - Classical HT Kingman: Comparison of exact and asymptotic results for $m_k(\alpha W)$. %, obtained via (\ref{6.7}) and numerical computation of the cumulants in (\ref{6.10}),
%and asymptotic result $k!(\frac12\,\sigma_{\alpha}^2)^k$.% with $\sigma_{\alpha}^2$ given in (\ref{6.9})
%, for the case that $\vart_U=5/2$, $\vart_V=1/2$ with $\alpha=1/10, 1/100, 1/1000$ and $k=1,2,3,4,5$.
}
\label{tbl:1}
\end{table}

\begin{figure}[ht]
\begin{tabular}{ccccc}
\includegraphics[width=0.31\textwidth]{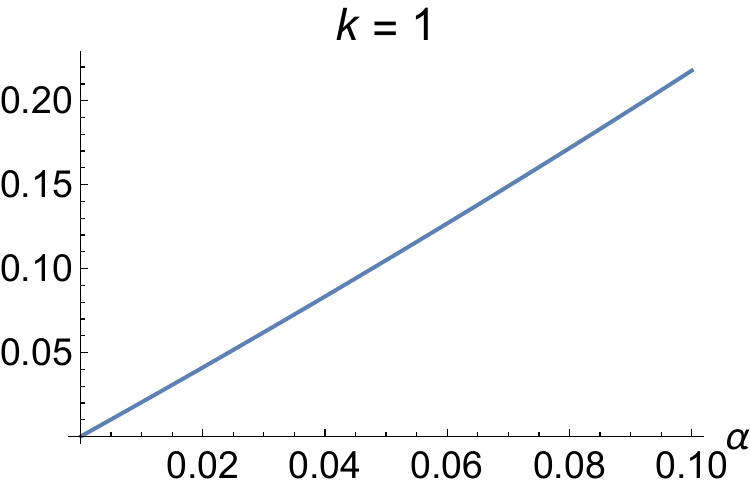} &
\includegraphics[width=0.31\textwidth]{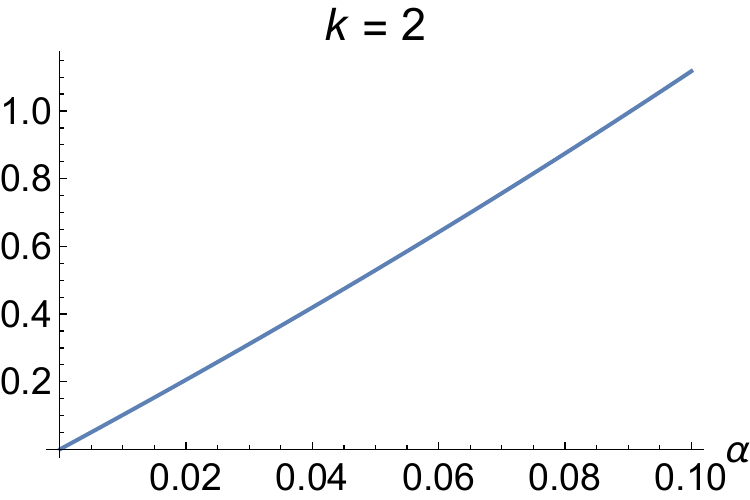} &
\includegraphics[width=0.31\textwidth]{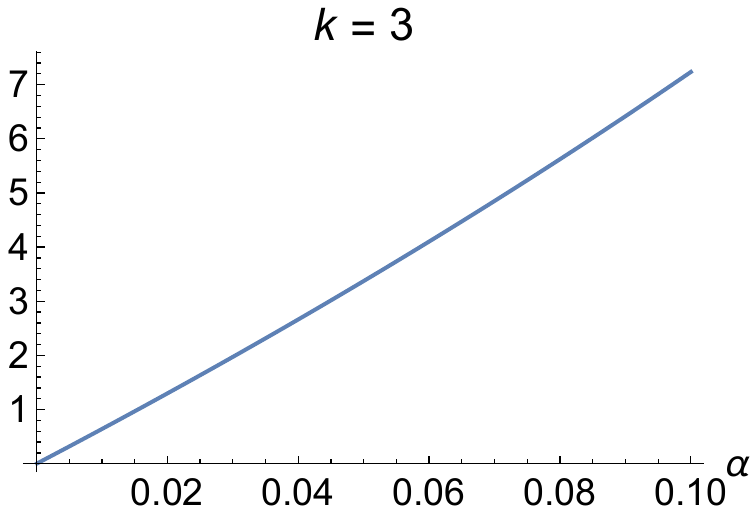}
\end{tabular}
\caption{Absolute error $k!(\frac12\,\sigma_{\alpha}^2)^k-m_k(\alpha W)$ for the classical Kingman heavy-traffic regime.}
\label{fig:error}
\end{figure}

We observe that the error %, the difference between the exact result and the asymptotic result,
behaves practically linearly with $\alpha$ for a fixed $k$. Thus, for this case, the error estimate $O(\alpha\,{\rm log}(1/\alpha))$ in Theorem~\ref{thm1} seems a factor ${\rm log}(1/\alpha)$ too pessimistic. To confirm this, we included a plot for the difference between the approximation (based on the asymptotic result) and the exact values for the first three moments. Indeed, Figure \ref{fig:error} confirms that the absolute error is almost completely linear in $\alpha$, meaning that the factor $\log(1/\alpha)$ is negligible here. Higher moments show the same type of behavior, but for reasons of compactness we have omitted the corresponding figures.

It is notable that for the first moment ($k=1$), the approximations systematically overestimate the exact values, which is known in the literature as Kingman's bound, see below \eqref{1.7a}. For higher moments, this is still the case for $\vart_U=5/2$, $\vart_V=1/2$, but choosing parameter values $\vart_U=1/2$, $\vart_V=5/2$ leads to underestimations. %As it turns out, this pattern is also observed in the nearly-deterministic regimes.
The $k$-behavior of the error for the three considered cases is markedly different, especially when $\alpha$ is not small.

\subsubsection{Nearly-deterministic Kingman regime} Take a fixed $\beta>0$, and let $1-\rho_n=\beta/n$. From Theorem~\ref{thm2}, we can approximate $\dE[W_n^k]$ by $k!\gamma_n^k$, with
\beq \label{6.12}
\gamma_n=(\sigma_V^2+\rho_n^{-2}\sigma_U^2)\,\rho_n/(2n(1-\rho_n))=(\vart_V+\rho_n^{-2}\vart_U)\,\rho_n/(2n(1-\rho_n)).
\eq
In Table~\ref{tbl:2}, we choose $\beta=1$. Note that the approximation, based on the asymptotic results, for $n=10,100,1000$ yields exactly the same numerical values as in the classical HT Kingman case with $\alpha=1/n$. This is due to the fact that accidentally $\frac12\sigma_{\alpha}^2=\gamma_n$ for our chosen parameter values.
\begin{table}[t]
\[
\begin{array}{|l|rr|rr|rr|rr|rr|}
\hline
\mc{11}{|c|}{%\text{Example 1  - Nearly-deterministic HT Kingman}~~ (
\vart_U=5/2, \vart_V=1/2}\\
\hline
& \mc{2}{c|}{n=10} & \mc{2}{c|}{n=100} & \mc{2}{c|}{n=1000}  & \mc{2}{c|}{n=10000} & \mc{2}{c|}{n=100000} \\
\hline
k & {\rm Exact} & {\rm Asymp} & {\rm Exact} & {\rm Asymp} & {\rm Exact} & {\rm Asymp}  & {\rm Exact} & {\rm Asymp} & {\rm Exact} & {\rm Asymp}\\
\hline
 1 & 1.212 & 1.614 & 1.404 & 1.510 & 1.469 & 1.501 & 1.490 & 1.500 & 1.497 & 1.500 \\
 2 & 3.572 & 5.209 & 4.205 & 4.561 & 4.405 & 4.506 & 4.470 & 4.501 & 4.490 & 4.500 \\
 3 & 15.711 & 25.222 & 18.880 & 20.663 & 19.819 & 20.291 & 20.114 & 20.254 & 20.207 & 20.250 \\
 4 & 92.084 & 162.819 & 113.030 & 124.814 & 118.888 & 121.825 & 120.680 & 121.532 & 121.241 & 121.503 \\
 5 & 674.652 & 1313.861 & 845.833 & 942.427 & 891.460 & 914.295 & 905.080 & 911.554 & 909.307 & 911.280 \\
\hline
\end{array}
\]
\[
\begin{array}{|l|rr|rr|rr|rr|rr|}
\hline
\mc{11}{|c|}{%\text{Example 1  - Nearly-deterministic HT Kingman}~~ (
\vart_U=1/2, \vart_V=5/2}\\
\hline
& \mc{2}{c|}{n=10} & \mc{2}{c|}{n=100} & \mc{2}{c|}{n=1000}  & \mc{2}{c|}{n=10000} & \mc{2}{c|}{n=100000} \\
\hline
k & {\rm Exact} & {\rm Asymp} & {\rm Exact} & {\rm Asymp} & {\rm Exact} & {\rm Asymp} & {\rm Exact} & {\rm Asymp} & {\rm Exact} & {\rm Asymp} \\
\hline
 1 & 1.151 & 1.403 & 1.397 & 1.490 & 1.468 & 1.499 & 1.490 & 1.500 & 1.497 & 1.500 \\
 2 & 3.551 & 3.936 & 4.204 & 4.440 & 4.405 & 4.494 & 4.470 & 4.499 & 4.490 & 4.500 \\
 3 & 16.371 & 16.562 & 18.962 & 19.849 & 19.828 & 20.210 & 20.115 & 20.246 & 20.207 & 20.250 \\
 4 & 100.564 & 92.932 & 114.027 & 118.300 & 118.993 & 121.176 & 120.691 & 121.468 & 121.242 & 121.497 \\
 5 & 772.098 & 651.817 & 857.114 & 881.352 & 892.646 & 908.217 & 905.200 & 910.946 & 909.320 & 911.220 \\
\hline
\end{array}
\]
\[
\begin{array}{|l|rr|rr|rr|rr|rr|}
\hline
\mc{11}{|c|}{%\text{Example 1  - Nearly-deterministic HT Kingman}~~ (
\vart_U=3/2, \vart_V=3/2}\\
\hline
& \mc{2}{c|}{n=10} & \mc{2}{c|}{n=100} & \mc{2}{c|}{n=1000}  & \mc{2}{c|}{n=10000} & \mc{2}{c|}{n=100000} \\
\hline
k & {\rm Exact} & {\rm Asymp} & {\rm Exact} & {\rm Asymp} & {\rm Exact} & {\rm Asymp} & {\rm Exact} & {\rm Asymp} & {\rm Exact} & {\rm Asymp} \\
\hline
 1 & 1.183 & 1.508 & 1.401 & 1.500 & 1.468 & 1.500 & 1.490 & 1.500 & 1.497 & 1.500 \\
 2 & 3.568 & 4.550 & 4.205 & 4.500 & 4.405 & 4.500 & 4.470 & 4.500 & 4.490 & 4.500 \\
 3 & 16.065 & 20.589 & 18.922 & 20.253 & 19.823 & 20.250 & 20.114 & 20.250 & 20.207 & 20.250 \\
 4 & 96.396 & 124.223 & 113.532 & 121.525 & 118.940 & 121.500 & 120.685 & 121.500 & 121.242 & 121.500 \\
 5 & 722.972 & 936.845 & 851.487 & 911.480 & 892.054 & 911.252 & 905.140 & 911.250 & 909.314 & 911.250 \\
\hline
\end{array}
\]
\caption{Example 1  - Nearly-deterministic HT Kingman: Comparison of exact and asymptotic results for $m_k(W_n)$.}
\label{tbl:2}
\end{table}

We observe that the error, the difference between the exact result and the asymptotic result, decays like $1/\sqrt{n}$ for a fixed $k$ (this decay behavior becomes more manifest when $n$ is further increased). This is in reasonable agreement with the error estimate $O({\rm log}\,n/\sqrt{n})$ that is given in Theorem~\ref{thm2}, (\ref{1.9}). The $k$-behavior of the error for the three considered cases is markedly different, especially for low $n$.

\subsubsection{Nearly-deterministic Gaussian regime}  We have now to invoke the whole machinery of the dedicated saddle point method. We have
\beq \label{6.14}
h(\zeta)={-}k_V\,{\rm log}(1+\vart_V\zeta)-k_U\,{\rm log}(1-\vart_U\zeta/\rho),~~~~~~\frac{-1}{\vart_V}<{\rm Re}(\zeta)<\frac{\rho}{\vart_U},
\eq
\beq \label{6.15}
\zeta_{sp}={-}\,\frac{1-\rho}{\vart_V+\vart_U},~~~~~~h''(\zeta_{sp})=\frac{(\vart_U+\vart_V)^3}{(\vart_U+\rho\vart_V)^2}, ~~~~~~ d_2={-}\frac{\vart_U^2-\vart_V^2}{3(\vart_U+\rho\vart_V)}.
\eq
Take a fixed $\beta>0$, and let $1-\rho_n=\beta/\sqrt{n}$. We approximate the moments of the scaled waiting times $\frac{\sqrt{n}}{\sigma_n}\,W_n$ by the moments of the maximum of the Gaussian random walk with drift $-\beta$. These $m_k(\beta)$ can be computed from the cumulants $c_l(M_{\beta})$ of $M_{\beta}$, using (\ref{78}) and \eqref{eqn87} or \eqref{eqn88}.

\begin{table}[t]
\[
\begin{array}{|l|rrr|rrr|rrr|}
\hline
\mc{10}{|c|}{%\text{Example 1  - Nearly-deterministic HT Gaussian}~~ (
\vart_U=5/2, \vart_V=1/2}\\
\hline
& \mc{3}{c|}{n=10} & \mc{3}{c|}{n=100} & \mc{3}{c|}{n=1000} \\
\hline
k & {\rm Exact} & {\rm Asymp\ 1} & {\rm Asymp\ 2} & {\rm Exact} & {\rm Asymp\ 1} & {\rm Asymp\ 2} & {\rm Exact} & {\rm Asymp\ 1} & {\rm Asymp\ 2} \\
\hline
% 1 & 0.367 & 0.375 & 0.378 & 0.402 & 0.402 & 0.403 & 0.410 & 0.410 & 0.411 \\
% 2 & 0.590 & 0.650 & 0.614 & 0.707 & 0.720 & 0.709 & 0.740 & 0.743 & 0.740 \\
% 3 & 1.356 & 1.629 & 1.376 & 1.794 & 1.870 & 1.796 & 1.928 & 1.951 & 1.929 \\
% 4 & 4.105 & 5.372 & 3.857 & 5.991 & 6.398 & 5.971 & 6.624 & 6.752 & 6.622 \\
% 5 & 15.483 & 22.057 & 12.501 & 24.931 & 27.267 & 24.664 & 28.346 & 29.106 & 28.320 \\
 1 & 0.3674 & 0.3748 & 0.3776 & 0.4021 & 0.4015 & 0.4030 & 0.4105 & 0.4100 & 0.4106 \\
 2 & 0.5898 & 0.6499 & 0.6136 & 0.7071 & 0.7202 & 0.7094 & 0.7396 & 0.7432 & 0.7398 \\
 3 & 1.3565 & 1.6288 & 1.3756 & 1.7936 & 1.8703 & 1.7963 & 1.9283 & 1.9514 & 1.9286 \\
 4 & 4.1049 & 5.3717 & 3.8574 & 5.9908 & 6.3985 & 5.9714 & 6.6235 & 6.7521 & 6.6218 \\
 5 & 15.4834 & 22.0568 & 12.5010 & 24.9315 & 27.2673 & 24.6639 & 28.3460 & 29.1063 & 28.3202 \\
 \hline
\end{array}
\]
\[
\begin{array}{|l|rrr|rrr|rrr|}
\hline
\mc{10}{|c|}{%\text{Example 1  - Nearly-deterministic HT Gaussian}~~ (
\vart_U=1/2, \vart_V=5/2}\\
\hline
& \mc{3}{c|}{n=10} & \mc{3}{c|}{n=100} & \mc{3}{c|}{n=1000} \\
\hline
k & {\rm Exact} & {\rm Asymp\ 1} & {\rm Asymp\ 2} & {\rm Exact} & {\rm Asymp\ 1} & {\rm Asymp\ 2} & {\rm Exact} & {\rm Asymp\ 1} & {\rm Asymp\ 2} \\
\hline
% 1 & 0.223 & 0.226 & 0.228 & 0.351 & 0.352 & 0.351 & 0.394 & 0.394 & 0.394 \\
% 2 & 0.324 & 0.314 & 0.337 & 0.601 & 0.593 & 0.603 & 0.704 & 0.701 & 0.704 \\
% 3 & 0.690 & 0.622 & 0.704 & 1.492 & 1.441 & 1.495 & 1.824 & 1.803 & 1.824 \\
% 4 & 1.936 & 1.608 & 1.892 & 4.880 & 4.602 & 4.868 & 6.224 & 6.109 & 6.222 \\
% 5 & 6.747 & 5.146 & 6.206 & 19.862 & 18.287 & 19.696 & 26.453 & 25.779 & 26.430 \\
 1 & 0.2227 & 0.2259 & 0.2282 & 0.3506 & 0.3523 & 0.3514 & 0.3938 & 0.3943 & 0.3938 \\
 2 & 0.3242 & 0.3139 & 0.3367 & 0.6006 & 0.5931 & 0.6025 & 0.7039 & 0.7009 & 0.7041 \\
 3 & 0.6903 & 0.6223 & 0.7043 & 1.4922 & 1.4410 & 1.4947 & 1.8236 & 1.8032 & 1.8239 \\
 4 & 1.9359 & 1.6078 & 1.8922 & 4.8802 & 4.6016 & 4.8685 & 6.2235 & 6.1091 & 6.2220 \\
 5 & 6.7472 & 5.1462 & 6.2056 & 19.8621 & 18.2869 & 19.6959 & 26.4527 & 25.7788 & 26.4304 \\
  \hline
\end{array}
\]

\[
\begin{array}{|l|rrr|rrr|rrr|}
\hline
\mc{10}{|c|}{%\text{Example 1  - Nearly-deterministic HT Gaussian}~~ (
\vart_U=3/2, \vart_V=3/2}\\
\hline
& \mc{3}{c|}{n=10} & \mc{3}{c|}{n=100} & \mc{3}{c|}{n=1000} \\
\hline
k & {\rm Exact} & {\rm Asymp\ 1} & {\rm Asymp\ 2} & {\rm Exact} & {\rm Asymp\ 1} & {\rm Asymp\ 2} & {\rm Exact} & {\rm Asymp\ 1} & {\rm Asymp\ 2} \\
\hline
% 1 & 0.292 & 0.299 & 0.299 & 0.376 & 0.377 & 0.377 & 0.402 & 0.402 & 0.402 \\
% 2 & 0.450 & 0.466 & 0.466 & 0.653 & 0.655 & 0.655 & 0.722 & 0.722 & 0.722 \\
% 3 & 1.006 & 1.045 & 1.045 & 1.641 & 1.646 & 1.646 & 1.876 & 1.876 & 1.876 \\
% 4 & 2.956 & 3.071 & 3.071 & 5.427 & 5.444 & 5.444 & 6.423 & 6.425 & 6.425 \\
% 5 & 10.798 & 11.217 & 11.217 & 22.348 & 22.418 & 22.418 & 27.394 & 27.402 & 27.402 \\
 1 & 0.2919 & 0.2985 & 0.2985 & 0.3761 & 0.3768 & 0.3768 & 0.4021 & 0.4022 & 0.4022 \\
 2 & 0.4504 & 0.4660 & 0.4660 & 0.6532 & 0.6550 & 0.6550 & 0.7217 & 0.7219 & 0.7219 \\
 3 & 1.0062 & 1.0449 & 1.0449 & 1.6410 & 1.6462 & 1.6462 & 1.8757 & 1.8763 & 1.8763 \\
 4 & 2.9556 & 3.0712 & 3.0712 & 5.4270 & 5.4442 & 5.4442 & 6.4226 & 6.4245 & 6.4245 \\
 5 & 10.7979 & 11.2169 & 11.2169 & 22.3476 & 22.4176 & 22.4176 & 27.3938 & 27.4020 & 27.4020 \\
\hline
\end{array}
\]
\caption{Example 1  - Nearly-deterministic HT Gaussian: Comparison of exact results for $m_k(\frac{\sqrt{n}}{\sigma_n}\, W_n)$ and the asymptotic results for the case that $\vart_U=5/2$, $\vart_V=1/2$ with $\beta=1$, $n=10, 100, 1000$ and $k=1, 2, 3, 4, 5$.  The two entries in the Asymp-columns give, for a particular $k$, the asymptotic result from (\ref{eqn86}) and (\ref{eqn89}) in that order.}
\label{tbl:3}
\end{table}

The results, shown in Table~\ref{tbl:3}, indicate that the absolute error behaves quite accurately as $O(1/\sqrt{n})$ and $O(1/n)$, respectively, for the two asymptotic estimates. It is also interesting to see that the two asymptotic estimates perform equally well for low $k$ and $n$, while the refined asymptotic estimate based on \eqref{eqn86} outperforms the asymptotic estimate based on \eqref{eqn89} in all other cases.

Notice that the two approximations yield the same results when $\vart_U=\vart_V=3/2$ in Table~\ref{tbl:3}. This is caused by $U$ and $V$ having the same Gamma distribution, resulting in $h'''(\zeta_{sp})=0$. It follows that $\phi_n=0$ and $B_n=\beta_n$ in Theorem~\ref{thm4}. Consequently, the refined approximation and the original approximation are equal here.

\section{Conclusions} \label{sec8}

We have presented, using Kingman's transform method based on Pollaczek's formula for the Laplace transform of the steady-state waiting time distribution of the ${\rm GI}/{\rm G}/1$ queue, several heavy-traffic limit theorems. Under the assumption that the distribution of both the service time and the interarrival time have a Laplace transform, analytic in an open strip containing the imaginary axis, these heavy-traffic limit results can be shown to be valid on the level of transforms in a full neighborhood of the origin, with error assessment. As a consequence, there is convergence for all moments, with a corresponding error assessment. We have considered the classical heavy-traffic regime in which the transform of the steady-state waiting time distribution converges, after appropriate scaling, to the transform of an exponentially distributed random variable as the system load $\rho=1-\alpha$ tends to 1 (Kingman-type result), with error shown to be bounded as $O(\alpha\,{\rm log}(1/\alpha))$ as $\alpha\downarrow0$. We also have considered nearly deterministic queues (obtained through cyclic thinning) in two different regimes, viz.\ where the system's load $\rho_n$ satisfies $1-\rho_n\asymp 1/n$ and $1-\rho_n\asymp 1/\sqrt{n}$, respectively, as the thinning factor $n$ tends to infinity. The first regime allows a result of the Kingman-type, viz.\ convergence in terms of transforms to an exponential distribution, with error bounded as $O(\frac{1}{\sqrt{n}}\,{\rm log}\,n)$. The second regime allows, after appropriate scaling, convergence in terms of transforms to the maximum of the Gaussian random walk with a specific negative drift. For this case, we have shown $O(1/\sqrt{n})$ error behavior of transforms and moments. The latter result is refined, so as to yield $O(1/n)$ errors, by judicious choice of the drift parameter, as well as by using a weakly non-linear transformation of the Laplace variable in the transform of the maximum of the Gaussian random walk. Our asymptotic results have been illustrated and compared with results obtained by numerical integration for various different combinations of service time and interarrival time distributions. For all regimes, the heavy-traffic limits for the moments of the waiting time distribution were also found to be good asymptotic approximations for the exact values.

\bibliographystyle{abbrv}
\bibliography{kingman}

\begin{thebibliography}{10}

\bibitem{ref1}
J.~Abate, G.~L. Choudhury, and W.~Whitt.
\newblock Calculation of the {GI/G/1} waiting time distribution and its
  cumulants from {Pollaczek's} formulas.
\newblock {\em Archiv f\"ur Elektrotechnik and \"Ubertragungstechnik},
  47(5/6):311--321, 1993.

\bibitem{abatewhitt92}
J.~Abate and W.~Whitt.
\newblock The {Fourier}-series method for inverting transforms of probability
  distributions.
\newblock {\em Queueing Systems}, 10:5--88, 1992.

\bibitem{ref2}
J.~Blanchet and P.~W. Glynn.
\newblock Complete corrected diffusion approximations for the maximum of a
  random walk.
\newblock {\em The Annals of Applied Probability}, 16(2):951--983, 2006.

\bibitem{boon2021optimal}
M.~Boon, G.~Janssen, J.~van Leeuwaarden, and R.~Timmerman.
\newblock Optimal capacity allocation for heavy-traffic fixed-cycle
  traffic-light queues and intersections.
\newblock {\em arXiv preprint arXiv:2104.04303}, 2021.

\bibitem{boon2019pollaczek}
M.~A. Boon, A.~Janssen, J.~S. van Leeuwaarden, and R.~W. Timmerman.
\newblock Pollaczek contour integrals for the fixed-cycle traffic-light queue.
\newblock {\em Queueing Systems}, 91(1):89--111, 2019.

\bibitem{ref3}
A.~A. Borovkov.
\newblock Some limit theorems in the theory of mass service.
\newblock {\em Theory of Probability and Its Applications}, 10(3):375--400,
  1965.

\bibitem{ref3b}
O.~J. Boxma and J.~W. Cohen.
\newblock Heavy-traffic analysis for the {GI/G/1} queue with heavy-tailed
  distributions.
\newblock {\em Queueing Systems}, 33:177--204, 1999.

\bibitem{ref4}
S.~L. Brumelle.
\newblock Some inequalities for parallel-server queues.
\newblock {\em Operations Research}, 19(2):402--413, 1971.

\bibitem{ref5}
J.~T. Chang and Y.~Peres.
\newblock Ladder heights, {Gaussian} random walks and the {Riemann} zeta
  function.
\newblock {\em Annals of Probability}, 25(2):787--802, 1997.

\bibitem{ref5a}
Y.~Chen and W.~Whitt.
\newblock Algorithms for the upper bound mean waiting time in the {GI/GI/1}
  queue.
\newblock {\em Queueing Systems}, 94:327–356, 2020.

\bibitem{cohen1993complex}
J.~Cohen.
\newblock Complex functions in queueing theory.
\newblock {\em AEU. Archiv f{\"u}r Elektronik und {\"U}bertragungstechnik},
  47(5-6):300--310, Pollaczek memorial volume, 1993.

\bibitem{cohen2012single}
J.~W. Cohen.
\newblock {\em The {S}ingle {S}erver {Q}ueue}.
\newblock North-Holland, Elsevier, 2nd edition, 1982.

\bibitem{ref6}
D.~L. Iglehart and W.~Whitt.
\newblock Multiple channel queues in heavy traffic.
\newblock {\em Advances in Applied Probability}, 2(2):355--369, 1970.

\bibitem{ref7}
A.~J. E.~M. Janssen and J.~S.~H. van Leeuwaarden.
\newblock On {Lerch's} transcendent and the {Gaussian} random walk.
\newblock {\em The Annals of Applied Probability}, 17(2):421--439, 2006.

\bibitem{ref8}
A.~J. E.~M. Janssen and J.~S.~H. van Leeuwaarden.
\newblock Cumulants of the maximum of the {Gaussian} random walk.
\newblock {\em Stochastic Processes and their Applications}, 117:1928--1959,
  2007.

\bibitem{mds}
A.~J. E.~M. Janssen and J.~S.~H. van Leeuwaarden.
\newblock Back to the roots of the {M/D/$s$} queue and the works of {Erlang},
  {Crommelin} and {Pollaczek}.
\newblock {\em Statistica Neerlandica}, 62(3):299--313, 2008.

\bibitem{ref9}
A.~J. E.~M. Janssen, J.~S.~H. van Leeuwaarden, and B.~W.~J. Mathijsen.
\newblock Novel heavy-traffic regimes for large-scale service systems.
\newblock {\em SIAM Journal of Applied Mathematics}, 78(2):787--812, 2015.

\bibitem{ref10}
J.~F.~C. Kingman.
\newblock The single server queue in heavy traffic.
\newblock {\em Mathematical Proceedings of the Cambridge Philosophical
  Society}, 57(4):902--904, 1961.

\bibitem{ref11}
J.~F.~C. Kingman.
\newblock On queues in heavy traffic.
\newblock {\em Journal of the Royal Statistical Society: Series B},
  24(2):383--392, 1962.

\bibitem{ref11a}
J.~F.~C. Kingman.
\newblock Some inequalities for the queue {GI/G/1}.
\newblock {\em Biometrika}, 49(3/4):315--324, 1962.

\bibitem{ref12}
J.~K\"ollerstr\"om.
\newblock Heavy traffic theory for queues with several servers.
\newblock {\em Journal of Applied Probability}, 11(3):544--552, 1974.

\bibitem{ref13}
J.~K\"ollerstr\"om.
\newblock Heavy traffic theory for queues with several severs.
\newblock {\em Journal of Applied Probability}, 16(2):393--401, 1979.

\bibitem{ref14}
B.~W.~J. Mathijsen, A.~J. E.~M. Janssen, J.~S.~H. van Leeuwaarden, and
  B.~Zwart.
\newblock Robust heavy-traffic approximations for service systems facing
  overdispersed demand.
\newblock {\em Queueing Systems}, 90:257--289, 2018.

\bibitem{ref15}
G.~R. Newell.
\newblock {\em Approximate Stochastic Behavior of $n$-Server Service Systems
  with Large $n$}.
\newblock Number~87 in Lecture Notes in Economics and Mathematical Systems.
  Springer-Verlag, 1973.

\bibitem{ref16}
D.~Siegmund.
\newblock Corrected diffusion approximations in certain random walk problems.
\newblock {\em Advances in Applied Probability}, 11(4):701--719, 1979.

\bibitem{ref17}
K.~Sigman and W.~Whitt.
\newblock Heavy-traffic limits for nearly deterministic queues.
\newblock {\em Journal of Applied Probability}, 48(3):657--678, 2011.

\bibitem{ref18}
K.~Sigman and W.~Whitt.
\newblock Heavy-traffic limits for nearly deterministic queues: Stationary
  distributions.
\newblock {\em Queueing Systems}, 69:145--173, 2011.

\bibitem{ref19}
W.~Whitt.
\newblock {\em Heavy Traffic Limit Theorems for Queues: {A} Survey}.
\newblock Number~98 in Lecture Notes in Economics and Mathematical Systems.
  Springer-Verlag, 19734.

\end{thebibliography}

\appendix

\section*{Appendix}

\section{Finishing the proof of Theorem~\ref{thm4}} \label{appA}

We approximate the front factor $FF$ in (\ref{5.17}) by using a linear approximation of $\zeta'(v)$ and a quadratic approximation of $\zeta(v)$ from (\ref{5.14}). Writing $\sigma=\sigma_n$ as we did in Section~\ref{sec5}, we have
\beq \label{A1}
\zeta'(y/\sigma\sqrt{n})=i-2d_2y/\sigma\sqrt{n}+O(y^2/n)~,
\eq
\beq \label{A2}
\sigma\sqrt{n}\,\zeta(y/\sigma\sqrt{n})={-}\beta+iy-\varp y^2+O(y^3/n)~,
\eq
where we have written $\beta=\beta_n$ and $\varp=\varp_n=d_2/\sigma_n\sqrt{n}$. Using this in (\ref{5.17}), we get
\begin{eqnarray} \label{A3}
FF & = & \frac{-s(1+2i\varp y+O(y^2/n))}{({-}\beta+iy-\varp y^2+O(y^3/n))(s+\beta-iy+\varp y^2+O(y^3/n))} \nonumber \\[3mm]
& = & {-}s\,\frac{1+i\varp y}{-\beta+iy-\varp y^2}~\frac{1+i\varp y}{s+\beta-iy+\varp y^2}\,(1+O(y^2/n))~.
\end{eqnarray}
Now
\begin{eqnarray} \label{A4}
\frac{-\beta+iy-\varp y^2}{1+i\varp y} & = & ({-}\beta+iy-\varp y^2)(1-i\varp y)+O(y^2/n) \nonumber \\[3mm]
& = & -\beta+iy(1+\beta\varp)+O(y^2/n)\nonumber \\[3mm]
& = & (1+\beta\varp)({-}B+iy)+O(y^2/n)~,
\end{eqnarray}
with $B=B_n=\beta_n/(1+\beta_n\varp_n)$, see (\ref{5.5}). Similarly,
\begin{eqnarray} \label{A5}
\hspace*{-9mm}\frac{s+\beta-iy+\varp y^2}{1+i\varp y} & \!= & \!s+\beta-iy(1+(s+\beta)\varp)+O(y^2/n) \nonumber \\[3mm]
& \!= & \!(1+(s+\beta)\varp)\,\Bigl(\frac{s+\beta}{1+(s+\beta)\varp}-iy\Bigr)+O(y^2/n)\,.
\end{eqnarray}
We further set, see (\ref{5.5}),
\beq \label{A6}
\frac{s+\beta}{1+(s+\beta)\varp}=R+\frac{\beta}{1+\beta\varp}=R+B~,
\eq
and we compute
\beq \label{A7}
R=\frac{s+\beta}{1+(s+\beta)\varp}-\frac{\beta}{1+\beta\varp}=\frac{s}{(1+(s+\beta)\varp)(1+\beta\varp)}~.
\eq
Therefore,
\begin{eqnarray} \label{A8}
FF & = & \frac{-s}{(1+\beta\varp)(1+(s+\beta)\varp)}~\frac{1+O(y^2/n)}{({-}B+iy)(R+B-iy)} \nonumber \\[3mm]
& = & \frac{R}{(B-iy)(R+B-iy)}~(1+O(y^2/n))~,
\end{eqnarray}
and this takes care of the factor in the integral in (\ref{5.16}).

We finally relate $n\,h(\zeta_{sp})$, occurring in the exponential in the integrand in (\ref{5.16}), and $B$. We have with $\beta=\beta_n$, $\sigma=\sigma_n$ given in (\ref{5.1}) and $d_2$ given by (\ref{5.3})
\begin{eqnarray} \label{A9}
B & = & \frac{\beta}{1+\beta\varp}=\frac{\beta_n}{1+\beta_n\,\tfrac{d_2}{\sigma\sqrt{n}}}=\frac{-\zeta_{sp}\,\sqrt{n\,h''(\zeta_{sp})}}{1+\tfrac{\zeta_{sp}\,h'''(\zeta_{sp})}{6h''(\zeta_{sp})}} \nonumber \\[3.5mm]
& = & -\zeta_{sp}\,\sqrt{n\,h''(\zeta_{sp})}\,\Bigl(1-\frac{\zeta_{sp}\,h'''(\zeta_{sp})}{6h''(\zeta_{sp})}+O\Bigl(\frac1n\Bigr)\Bigr)~.
\end{eqnarray}
On the other hand, from $h'(\zeta_{sp})=0=h(0)$ and
\beq \label{A10}
h(0)=h(\zeta_{sp})-\zeta_{sp}\,h'(\zeta_{sp})+\tfrac12\,\zeta_{sp}^2\,h''(\zeta_{sp})-\tfrac16\,\zeta_{sp}^3\,h'''(\zeta_{sp})+O\Bigl(\frac{1}{n^2}\Bigr)~,
\eq
we get
\begin{eqnarray} \label{A11}
-n\,h(\zeta_{sp}) & = & \tfrac12\,n\,\zeta_{sp}^2\,h''(\zeta_{sp})-\tfrac16\,n\,\zeta_{sp}^3\,h'''(\zeta_{sp})+O\Bigl(\frac1n\Bigr) \nonumber \\[3mm]
& = & \tfrac12\,n\,\zeta_{sp}^2\,h''(\zeta_{sp})\Bigl(1-\frac{\zeta_{sp}\,h'''(\zeta_{sp})}{3h''(\zeta_{sp})}\Bigr)+O\Bigl(\frac1n\Bigr) \nonumber \\[3mm]
& = & \tfrac12\,B^2+O\Bigl(\frac1n\Bigr)~.
\end{eqnarray}
Using (\ref{A8}) and (\ref{A11}), valid uniformly in any bounded set of $s$ with ${\rm Re}(s)\geq{-}\frac12\,\beta_n$, in (\ref{5.16}), we get
\begin{eqnarray} \label{A12}
~& \mbox{} & {\rm log}(\dE\,\Bigl[\exp\Bigl({-}\,\frac{s\sqrt{n}}{\sigma_n}\,W_n\Bigr)\Bigr]) \nonumber \\[3mm]
& & =~\frac{1}{2\pi i}\,\il_{-R}^R\,\frac{R_n\,{\rm log}(1-e^{-\frac12 B_n^2-\frac12 y^2})}{(B_n-iy)(R_n+B_n-iy)}\,dy\,\Bigl(1+O\Bigl(\frac1n\Bigr)\Bigr)~,
\end{eqnarray}
where we have restored the $n$ in $\sigma_n$, $R_n$ and $B_n$. Then the proof of Theorem~\ref{thm4} can be finished in the same way as the proof of Theorem~\ref{thm3} was finished.

\section{Proof of Theorem \ref{thm5}} \label{appB}

The first line of (\ref{5.6}) is an immediate consequence of Theorem~\ref{thm4}, and can be rewritten as
\beq \label{B1}
\dE\,\Bigl[\Bigl(\frac{\sqrt{n}}{\sigma_n}\,W_n\Bigr)^k\Bigr]=\Bigl(\frac{d}{ds}\Bigr)^k\,[F(T_n(s))]_{s=0}+O\Bigl(\frac1n\Bigr)~,
\eq
where $F(s)=\dE\,[\exp(s\,M_{B_n})]$ and $T_n(s)={-}R_n({-}s)$ with $B_n$ and $R_n$ given in (\ref{5.5}).
We have from (\ref{5.5}), deleting the index $n$ from $\beta_n$ and $\varp_n$,
\beq \label{B2}
T_n(s)=\frac{s}{1+\beta\varp}~\frac{1}{1+\beta\varp-\varp s}=\sum_{r=1}^{\infty}\,\frac{\varp^{r-1}}{(1+\beta\varp)^{r+1}}\,s^r~.
\eq
Therefore, $T_n(0)=0$ and
\beq \label{B3}
T_n^{(r)}(0)=\frac{r!\,\varp^{r-1}}{(1+\beta\varp)^{r+1}}~,~~~~~~r=1,2,...~,
\eq
so that, in particular,
\beq \label{B4}
T_n'(0)=\frac{1}{(1+\beta\varp)^2}~,~~~~~~T_n''(0)=\frac{2\varp}{(1+\beta\varp)^3}~,
\eq
\beq \label{B5}
T_n^{(r)}(0)=O(\varp^{r-1})=O(n^{-\frac12(r-1)})~,
\eq
see below (\ref{5.5}).

We compute
\beq \label{B6}
V_k=\Bigl(\frac{d}{ds}\Bigr)^k\,[F(T_n(s))](0)~,
\eq
noting that $F^{(k)}(0)=m_k(B_n)$ for $k=1,2,...\,$, within absolute error $O(1/n)$. We have, writing $B=B_n$,
\beq \label{B7}
V_1=F'(T_n(s))\,T_n'(s)\Bigl|_{s=0}=F'(0)\,T_n'(0)=\frac{m_1(B)}{(1+\beta\varp)^2}~,
\eq
\begin{eqnarray} \label{B8}
V_2 & = & \frac{d}{ds}\,(F'(T_n(s))\,T_n'(s))\Bigl|_{s=0} \nonumber \\[3mm]
& = & (F''(T_n(s))(T_n'(s))^2+F'(T_n(s))T_n''(s))\Bigl|_{s=0} \nonumber \\[3mm]
& = & F''(0)(T_n'(0))^2+F'(0)\,T_n''(0)=\frac{m_2(B)}{(1+\beta\varp)^4}+\frac{2\varp\,m_1(B)}{(1+\beta\varp)^3}~,
\end{eqnarray}
and similarly
\begin{eqnarray} \label{B9}
V_3 & = & F'''(0)(T_n'(0))^3+3F''(0)\,T_n'(0)\,T_n''(0)+F'(0)\,T_n'''(0) \nonumber \\[3mm]
& = & \frac{m_3(B)}{(1+\beta\varp)^6}+\frac{6\varp\,m_2(B)}{(1+\beta\varp)^5}+O\Bigl(\frac1n\Bigr)~,
\end{eqnarray}
where (\ref{B4}--\ref{B5}) has been used. In general, one finds inductively
\begin{eqnarray} \label{B10}
V_k & = & (F^{(k)}(T_n(s))(T_n'(s))^k+\tfrac12\,k(k-1)\,F^{(k-1)}(T_n(s))(T_n'(s))^{k-1}T_n''(s)+...)\Bigl|_{s=0} \nonumber \\[3mm]
& = & \frac{m_k(B)}{(1+\beta\varp)^{2k}}+\frac{k(k-1)\,\varp\,m_{k-1}(B)}{(1+\beta\varp)^{2k-1}}+O\Bigl(\frac1n\Bigr)~,
\end{eqnarray}
and this gives the expression on the second line of (\ref{5.6}) after restoring the $n$ in $B$, $\beta$ and $\varp$.

We finally show that Theorem~\ref{thm3} gives qualitatively the same accurate estimates of the moments as Theorem~\ref{thm5} does when the third cumulants of $V$ and $U$ are equal. We have, abbreviating ``$k^\text{th}$ cumulant of'' by ``$c_k$'',
\begin{align}
h(\zeta)\,&=\,\log\left(\dE[e^{-\zeta V}]\right)+\log\left(\dE[e^{\zeta U/\rho}]\right)\nonumber\\
&=\,\sum_{k=1}^\infty \frac{\rho^{-k}c_k(U)+(-1)^{k}c_k(V)}{k!}\zeta^k.\label{eqn:B11}
\end{align}
With $1-\rho=1-\rho_n\asymp1/\sqrt{n}$, we have that $\zeta_{sp}=O(1/\sqrt{n})$. Assume that $c_3(V)=c_3(U)$. Then from \eqref{eqn:B11}
\begin{equation}\label{eqn:B12}
h'''(\zeta_{sp})=O(\rho^{-3}c_3(U)-c_3(V))+O(\zeta_{sp})=O(1/\sqrt{n}).
\end{equation}
Therefore, $\phi_n$ in Theorem~\ref{thm4} is $O(1/n)$, and $B_n=\beta_n+O(1/n)$. Hence, the leading term in \eqref{5.6} and the leading term in \eqref{1.10} agree with one another within an error $O(1/n)$, and so the first order term in \eqref{5.6} approximates the scaled moment at the left-hand side of \eqref{5.6} within an error $O(1/n)$. We conclude that the leading term in \eqref{1.10} approximates the scaled moment at the left-hand side of \eqref{5.6} and \eqref{1.10} within an error $O(1/n)$ as well.

%\end{document}

\section{Further numerical examples}\label{sec7}

In the main text of this paper we have provided one numerical example (referred to as Example~1). We now present more
numerical examples to illustrate the applicability and accuracy of the proposed approximations. Although the main focus of this paper is on the \emph{moments} of the waiting-time distribution, we also include a few results for the \emph{cumulative distribution functions}. We use the results of Section~\ref{sec6} to compute the moments in the classical HT Kingman case, the nearly deterministic HT Kingman case and the nearly deterministic HT Gaussian case.

\subsection{Example 2: $U$ Gamma distributed and $V$ Bates distributed}

In this example, we take a different distribution for the generic service times $V$, namely the mean of $m$ independent uniform $[1-\delta, 1+\delta]$ random variables, with pdf
proportional to
\begin{equation}
\sum _{k=0}^m (-1)^k \binom{m}{k} \left(\frac{x-1+\delta}{2 \delta }-\frac{k}{m}\right)^{m-1} \text{sgn}\left(\frac{x-1+\delta}{2 \delta }-\frac{k}{m}\right), \qquad 1-\delta\leq x\leq 1+\delta.
\end{equation}
This distribution, also known as the Bates distribution, has mean 1 and variance $\frac{\delta^2}{3m}$. We choose $m=4$ and $\delta=1$.
The unscaled interarrival times are, again, Gamma distributed with parameters $k_U=\frac25$ and $\vart_U=\frac52$. This gives
\begin{equation}
\psi(-\zeta)=\dE\,[e^{-\zeta(V-\frac{1}{\rho}U)}]= \left(\frac{e^{-\frac{(1+\delta) \zeta }{m}}-e^{-\frac{(1-\delta ) \zeta }{m}}}{-2\delta  \zeta/m
   }\right)^m (1-{\vart_U \zeta  }/{\rho })^{-k_U}
\end{equation}
for $-\infty<\text{Re}(\zeta)<\rho/\vart_U$. As in the previous example, Assumption~\ref{ass:delta} of Section~\ref{sec5} is easily checked.

\noindent\textbf{Classical HT Kingman.} Again, we approximate the moments of $\alpha W$ by \eqref{eqn80}, with
\beq \label{6.9b}
\sigma_{\alpha}^2=(\sigma_V^2+\rho^{-2}\sigma_U^2)\,\rho=\frac{\delta^2\rho}{3m}+\frac{k_U\vart_U^2}{\rho}=
\frac{\rho}{12}+\frac{5}{2\rho}~,
\eq
and $\rho=1-\alpha$, where we take $\alpha=1/10,1/100,1/1000$. The results, for the first five moments of $W$, are presented in Table~\ref{tbl:1-1b}.

\begin{table}[t]
\[
\begin{array}{|l|rr|rr|rr|}
%\hline
%\mc{7}{c}{\text{Example 2  - Classical HT Kingman}~~ (\vart_U=5/2, m=4, \delta=1)}\\
\hline
& \mc{2}{c|}{\alpha=1/10} & \mc{2}{c|}{\alpha=1/100} & \mc{2}{c|}{\alpha=1/1000} \\
\hline
k & {\rm Exact} & {\rm Asymp} & {\rm Exact} & {\rm Asymp} & {\rm Exact} & {\rm Asymp} \\
\hline
 1 & 1.205 & 1.426 & 1.283 & 1.304 & 1.291 & 1.293 \\
 2 & 3.014 & 4.069 & 3.304 & 3.400 & 3.334 & 3.343 \\
 3 & 11.300 & 17.413 & 12.764 & 13.300 & 12.914 & 12.966 \\
 4 & 56.488 & 99.349 & 65.742 & 69.368 & 66.699 & 67.056 \\
 5 & 352.976 & 708.551 & 423.246 & 452.234 & 430.629 & 433.477 \\
\hline
\end{array}
\]
\caption{Example 2  - Classical HT Kingman: Comparison of exact result for $m_k(\alpha W)$ %, obtained via (\ref{6.7}) and numerical computation of the cumulants in (\ref{6.10}),
and asymptotic result $k!(\frac12\,\sigma_{\alpha}^2)^k$% with $\sigma_{\alpha}^2$ given in (\ref{6.9})
%, for the case that $\vart_U=5/2$, $m=4$, $\delta = 1$
 with $\alpha=1/10,1/100,1/1000$ and $k=1,2,3,4,5$.}
\label{tbl:1-1b}
\end{table}

\noindent\textbf{Nearly-deterministic HT Kingman.} With $\beta=1$ and $\rho_n=1-\beta/n$, we get
\beq \label{6.12b}
\gamma_n=(\sigma_V^2+\rho_n^{-2}\sigma_U^2)\,\rho_n/(2n(1-\rho_n))=
\frac{1}{2n(1-\rho_n)}\left(\frac{\delta^2\rho_n}{3m}+\frac{k_U\vart_U^2}{\rho_n}\right)=
\frac{\rho}{24}+\frac{5}{4\rho}~.
\eq
In Table \ref{tbl:2-1b}, we show results for for $n=10, 100, 1000, 10000, 100000$, which corresponds to $\rho_n=0.9, 0.99$, $0.999$, $0.9999$, $0.99999$, respectively. As in Example 1, we note that the approximation yields the exact same numerical values as in the classical HT Kingman case with $\alpha=1/n$. This is again due to the fact that $\frac12\sigma_{\alpha}^2=\gamma_n$ for $\beta=1$.

\begin{table}[t]
\[
\begin{array}{|l|rr|rr|rr|rr|rr|}
%\hline
%\mc{7}{c}{\text{Example 2  - Nearly-deterministic HT Kingman}~~ (\vart_U=1/2, m=4, \delta=1)}\\
\hline
& \mc{2}{c|}{n=10} & \mc{2}{c|}{n=100} & \mc{2}{c|}{n=1000}& \mc{2}{c|}{n=10000} & \mc{2}{c|}{n=100000} \\
\hline
k & {\rm Exact} & {\rm Asymp} & {\rm Exact} & {\rm Asymp} & {\rm Exact} & {\rm Asymp} & {\rm Exact} & {\rm Asymp} & {\rm Exact} & {\rm Asymp} \\
\hline
 1 & 1.035 & 1.426 & 1.204 & 1.304 & 1.263 & 1.293 & 1.282 & 1.292 & 1.289 & 1.292 \\
 2 & 2.604 & 4.069 & 3.102 & 3.400 & 3.261 & 3.343 & 3.313 & 3.337 & 3.329 & 3.337 \\
 3 & 9.767 & 17.413 & 11.982 & 13.300 & 12.633 & 12.966 & 12.836 & 12.934 & 12.900 & 12.930 \\
 4 & 48.827 & 99.349 & 61.710 & 69.368 & 65.249 & 67.056 & 66.319 & 66.831 & 66.652 & 66.808 \\
 5 & 305.107 & 708.551 & 397.293 & 452.234 & 421.271 & 433.477 & 428.296 & 431.655 & 430.461 & 431.473 \\
\hline
\end{array}
\]
\caption{Example 2  - Nearly-deterministic HT Kingman: Comparison of exact result for $m_k(\alpha W)$ %, obtained via (\ref{6.7}) and numerical computation of the cumulants in (\ref{6.10}),
and asymptotic result $k!(\frac12\,\sigma_{\alpha}^2)^k$% with $\sigma_{\alpha}^2$ given in (\ref{6.9})
%, for the case that $\vart_U=5/2$, $m=4$, $\delta = 1$
 with $n=10,100,1000,10000,100000$ and $k=1,2,3,4,5$.}
\label{tbl:2-1b}
\end{table}

\noindent\textbf{Nearly-deterministic HT Gaussian.}
It is readily seen that
\begin{equation}
h(\zeta)=m\left(\log\frac{m}{2s}-\frac{\zeta}{m}+\log\left(\frac{e^{{-\delta \zeta }/{m}}-e^{{\delta  \zeta }/{m}}}{-\zeta}\right)\right)-\frac{1}{\vart_U}\log(1-\vart_U \zeta/\rho).
\end{equation}
Although it is straightforward to compute the derivatives of $h(\zeta)$, in this case there is no explicit expression for the saddle point $\zeta_{sp}$. It can be shown, though, that
\begin{equation}\label{zetaspapprox}
\zeta_{sp}=\frac{1/\rho-1}{\delta^2/3m+\vart/\rho^2}+O\big((1-\rho)^2\big).
\end{equation}
As a consequence, $\zeta_{sp}$ can be determined numerically, for example using a Newton iteration with the leading term on the right hand side of \eqref{zetaspapprox} as a starting point.
The exact and approximate values for the first five moments of the scaled waiting times are shown in Table \ref{tbl:3-1b}. In this case, $c_3(V)=0\neq 2\vart_U^2=c_3(U)$, and so the refined asymptotic result (Asymp 2) should outperform the standard asymptotic result (Asymp 1) at least for large $n$. This is confirmed by Table~\ref{tbl:3-1b}.
\begin{table}[t]
\[
\begin{array}{|l|rrr|rrr|rrr|}
%\hline
%\mc{10}{c}{\text{Example 2  - Nearly-deterministic HT Gaussian}~~ (\vart_U=1/2,  m=4, \delta=1)}\\
\hline
& \mc{3}{c|}{n=10} & \mc{3}{c|}{n=100} & \mc{3}{c|}{n=1000} \\
\hline
k & {\rm Exact} & {\rm Asymp\ 1} & {\rm Asymp\ 2} & {\rm Exact} & {\rm Asymp\ 1} & {\rm Asymp\ 2} & {\rm Exact} & {\rm Asymp\ 1} & {\rm Asymp\ 2} \\
\hline
% 1 & 0.345 & 0.354 & 0.355 & 0.359 & 0.358 & 0.359 & 0.360 & 0.360 & 0.360 \\
% 2 & 0.526 & 0.597 & 0.548 & 0.593 & 0.608 & 0.595 & 0.608 & 0.612 & 0.608 \\
% 3 & 1.134 & 1.452 & 1.120 & 1.405 & 1.490 & 1.405 & 1.477 & 1.503 & 1.477 \\
% 4 & 3.215 & 4.647 & 2.654 & 4.377 & 4.802 & 4.334 & 4.721 & 4.852 & 4.717 \\
% 5 & 11.360 & 18.507 & 5.875 & 16.973 & 19.254 & 16.543 & 18.779 & 19.498 & 18.739 \\
 1 & 0.3453 & 0.3537 & 0.3554 & 0.3586 & 0.3584 & 0.3595 & 0.3602 & 0.3599 & 0.3603 \\
 2 & 0.5255 & 0.5966 & 0.5483 & 0.5926 & 0.6083 & 0.5947 & 0.6077 & 0.6120 & 0.6079 \\
 3 & 1.1345 & 1.4523 & 1.1204 & 1.4051 & 1.4904 & 1.4048 & 1.4773 & 1.5027 & 1.4773 \\
 4 & 3.2152 & 4.6472 & 2.6535 & 4.3767 & 4.8017 & 4.3336 & 4.7210 & 4.8517 & 4.7171 \\
 5 & 11.3597 & 18.5069 & 5.8753 & 16.9729 & 19.2541 & 16.5425 & 18.7795 & 19.4975 & 18.7396 \\
\hline
\end{array}
\]
\caption{Example 2  - Nearly-deterministic HT Gaussian: Comparison of exact result for $m_k(\frac{\sqrt{n}}{\sigma_n}\, W_n)$ and the asymptotic results for %the case% that $\vart_U=5/2$, $m=4$, $\delta = \frac12$
% with $\beta=1$,
 $n=10,100,1000$ and $k=1,2,3,4,5$.  The two entries in the Asymp-columns give, for a particular $k$, the asymptotic result from (\ref{eqn86}) and (\ref{eqn89}) in that order.}
\label{tbl:3-1b}
\end{table}

\subsection{Example 3: $U$ Gamma distributed and $V$ lattice distributed}

In this example, the generic service times $V$ follow an $m$-point lattice distribution centered around 1 in the interval $[1-\delta, 1+\delta]$, with $0<\delta\leq 1$ and $m>1$. To be precise, $V$ takes values $(1-\delta) + k \frac{2\delta}{m-1}$ for $k=0,1,\dots,m-1$, each with probability $1/m$. The variance of $V$ is equal to $\frac{\delta^2(m+1)}{3(m-1)}$, which is always less than one, except when $m=2$ and $\delta=1$. In this example, we take $m=4$, $\delta = 1$, and $\vart_U=\frac52$. Note that this means that service times are equal to zero with probability $1/4$. The interesting feature of this example is that, with $V$ having a lattice distribution, the conditions required for the Nearly-deterministic Gaussian regime are more subtle. Nevertheless, we still get accurate approximations when applying our saddle point technique, as shown in Table~\ref{tbl:3-1c}. For completeness, we also present the results for the classical Kingman-type heavy traffic regime and for the nearly-deterministic Kingman-type heavy traffic regime.

\noindent\textbf{Classical HT Kingman.} Again, we approximate the moments of $\alpha W$ by \eqref{eqn80}, with
\beq \label{6.9c}
\sigma_{\alpha}^2=(\sigma_V^2+\rho^{-2}\sigma_U^2)\,\rho=\frac{\delta^2(m+1)\rho}{3(m-1)}+\frac{k_U\vart_U^2}{\rho}=
\frac{5\rho}{9}+\frac{5}{2\rho}~,
\eq
and $\rho=1-\alpha$, where we take $\alpha=1/10,1/100,1/1000$. The results, for the first five moments of $W$, are presented in Table~\ref{tbl:1-1c}.

\begin{table}[t]
\[
\begin{array}{|l|rr|rr|rr|}
%\hline
%\mc{7}{c}{\text{Example 3  - Classical HT Kingman}~~ (\vart_U=5/2, m=4, \delta=1)}\\
\hline
& \mc{2}{c|}{\alpha=1/10} & \mc{2}{c|}{\alpha=1/100} & \mc{2}{c|}{\alpha=1/1000} \\
\hline
k & {\rm Exact} & {\rm Asymp} & {\rm Exact} & {\rm Asymp} & {\rm Exact} & {\rm Asymp} \\
\hline
 1 & 1.423 & 1.639 & 1.517 & 1.538 & 1.527 & 1.529 \\
 2 & 4.234 & 5.372 & 4.625 & 4.729 & 4.664 & 4.674 \\
 3 & 18.888 & 26.412 & 21.142 & 21.812 & 21.371 & 21.437 \\
 4 & 112.340 & 173.145 & 128.868 & 134.157 & 130.564 & 131.087 \\
 5 & 835.213 & 1418.825 & 981.875 & 1031.419 & 997.109 & 1001.995 \\
\hline
\end{array}
\]
\caption{Example 3  - Classical HT Kingman: Comparison of exact result for $m_k(\alpha W)$ %, obtained via (\ref{6.7}) and numerical computation of the cumulants in (\ref{6.10}),
and asymptotic result $k!(\frac12\,\sigma_{\alpha}^2)^k$% with $\sigma_{\alpha}^2$ given in (\ref{6.9})
%, for the case that %$\vart_U=5/2$, $m=4$, $\delta = 1$
with $\alpha=1/10,1/100,1/1000$ and $k=1,2,3,4,5$.}
\label{tbl:1-1c}
\end{table}

\noindent\textbf{Nearly-deterministic HT Kingman.} With $\beta=1$ and $\rho_n=1-\beta/n$, we get
\beq \label{6.12b}
\gamma_n=(\sigma_V^2+\rho_n^{-2}\sigma_U^2)\,\rho_n/(2n(1-\rho_n))=\frac{5\rho}{9}+\frac{5}{2\rho}~.
\eq
As in Example 2, we see that $\gamma_n$ in the Nearly-deterministic HT Kingman regime is equal to $\frac12\sigma_{\alpha}^2$, due to the fact that $\beta=1$, resulting in the same limiting values (shown in Table~\ref{tbl:2-1c}).

\begin{table}[t]
\[
\begin{array}{|l|rr|rr|rr|rr|rr|}
%\hline
%\mc{7}{c}{\text{Example 3  - Nearly-deterministic HT Kingman}~~ (\vart_U=5/2, m=4, \delta=1)}\\
\hline
& \mc{2}{c|}{n=10} & \mc{2}{c|}{n=100} & \mc{2}{c|}{n=1000} & \mc{2}{c|}{n=10000}  & \mc{2}{c|}{n=100000} \\
\hline
k & {\rm Exact} & {\rm Asymp} & {\rm Exact} & {\rm Asymp} & {\rm Exact} & {\rm Asymp} & {\rm Exact} & {\rm Asymp}  & {\rm Exact} & {\rm Asymp} \\
\hline
 1 & 1.235 & 1.639 & 1.430 & 1.538 & 1.496 & 1.529 & 1.518 & 1.528 & 1.525 & 1.528 \\
 2 & 3.692 & 5.372 & 4.362 & 4.729 & 4.570 & 4.674 & 4.637 & 4.669 & 4.658 & 4.668 \\
 3 & 16.476 & 26.412 & 19.942 & 21.812 & 20.943 & 21.437 & 21.253 & 21.400 & 21.351 & 21.396 \\
 4 & 97.997 & 173.145 & 121.552 & 134.157 & 127.950 & 131.087 & 129.877 & 130.786 & 130.477 & 130.756 \\
 5 & 728.576 & 1418.825 & 926.135 & 1031.419 & 977.143 & 1001.995 & 992.089 & 999.126 & 996.697 & 998.840 \\
\hline
\end{array}
\]
\caption{Example 3  - Nearly-deterministic HT Kingman: Comparison of exact result for $m_k(W_n)$ %, obtained via (\ref{6.7}) and numerical computation of the cumulants in (\ref{6.10}),
and asymptotic result $k!\gamma_n^k$% with $\sigma_{\alpha}^2$ given in (\ref{6.9})
%, for the case that $\vart_U=5/2$, $m=4$, $\delta = 1$
 with $n=10, 100, 1000, 10000, 100000$ and $k=1, 2, 3, 4, 5$.}
\label{tbl:2-1c}
\end{table}

\noindent\textbf{Nearly-deterministic HT Gaussian.}
To use the saddle point method, we first compute $h(\zeta)$:
\begin{align}
h(\zeta)%=\log \psi(-\zeta)
&= \log\left(\frac1m\sum_{k=0}^{m-1} e^{-\zeta\big(1-\delta + k \frac{2\delta}{m-1}\big)}\right)-\frac{1}{\vart_U}\log(1-\vart_U \zeta/\rho)\nonumber\\
&= -\zeta(1-\delta)+ \log\left(\frac1m\sum_{k=0}^{m-1} e^{-\zeta k \frac{2\delta}{m-1}}\right)-\frac{1}{\vart_U}\log(1-\vart_U \zeta/\rho).
\end{align}
As in the previous example, there is no closed-from expression for $\zeta_{sp}$, but with all the derivatives of $h(\zeta)$ being analytic, numerical methods can easily find the saddle point of $h(\zeta)$. In fact, our numerical methods encountered no issues computing the exact values and the two approximations for the first five moments of $\frac{\sqrt{n}}{\sigma_n}W_n$. The results are shown in Table~\ref{tbl:3-1c}. In this example, with $c_3(V)=0$ and $c_3(U)=2\vart_U^2>0$, the refined approximation gives significantly better results than the standard approximation.
\begin{table}[t]
\[
\begin{array}{|l|rrr|rrr|rrr|}
%\hline
%\mc{10}{c}{\text{Example 3  - Nearly-deterministic HT Gaussian}~~ (\vart_U=5/2,  m=4, \delta=1)}\\
\hline
& \mc{3}{c|}{n=10} & \mc{3}{c|}{n=100} & \mc{3}{c|}{n=1000} \\
\hline
k & {\rm Exact} & {\rm Asymp\ 1} & {\rm Asymp\ 2} & {\rm Exact} & {\rm Asymp\ 1} & {\rm Asymp\ 2} & {\rm Exact} & {\rm Asymp\ 1} & {\rm Asymp\ 2} \\
\hline
% 1 & 0.373 & 0.380 & 0.383 & 0.408 & 0.407 & 0.409 & 0.417 & 0.417 & 0.417 \\
% 2 & 0.602 & 0.662 & 0.624 & 0.723 & 0.736 & 0.725 & 0.757 & 0.761 & 0.758 \\
% 3 & 1.389 & 1.671 & 1.398 & 1.846 & 1.926 & 1.848 & 1.991 & 2.015 & 1.991 \\
% 4 & 4.216 & 5.546 & 3.897 & 6.213 & 6.639 & 6.188 & 6.898 & 7.032 & 6.896 \\
% 5 & 15.944 & 22.928 & 12.401 & 26.049 & 28.514 & 25.744 & 29.776 & 30.578 & 29.746 \\
 1 & 0.3734 & 0.3796 & 0.3827 & 0.4082 & 0.4074 & 0.4090 & 0.4171 & 0.4166 & 0.4171 \\
 2 & 0.6023 & 0.6622 & 0.6237 & 0.7227 & 0.7358 & 0.7248 & 0.7574 & 0.7610 & 0.7576 \\
 3 & 1.3894 & 1.6705 & 1.3979 & 1.8464 & 1.9255 & 1.8481 & 1.9911 & 2.0150 & 1.9913 \\
 4 & 4.2156 & 5.5461 & 3.8966 & 6.2127 & 6.6388 & 6.1882 & 6.8978 & 7.0322 & 6.8955 \\
 5 & 15.9440 & 22.9275 & 12.4014 & 26.0485 & 28.5145 & 25.7437 & 29.7756 & 30.5780 & 29.7463 \\
\hline
\end{array}
\]
\caption{Example 3  - Nearly-deterministic HT Gaussian: Comparison of exact result for $m_k(\frac{\sqrt{n}}{\sigma_n}\, W_n)$ and the asymptotic results% for the case that $\vart_U=5/2$, $m=4$, $\delta = 1$ with $\beta=1$,
with $n=10, 100, 1000$ and $k=1, 2, 3, 4, 5$.  The two entries in the Asymp-columns give, for a particular $k$, the asymptotic result from (\ref{eqn86}) and (\ref{eqn89}) in that order.}
\label{tbl:3-1c}
\end{table}

\subsection{Example 4: $U$ and $V$ lattice distributed}

In this example, \emph{both} the generic interarrival times $U$ and service times $V$ follow a two-point lattice distribution with mean 1:
\begin{equation}
U=\begin{cases}
a_U & \quad \text{ w.p. } p_U,\\
b_U & \quad \text{ w.p. } 1-p_U,
\end{cases}
\qquad V=\begin{cases}
a_V & \quad \text{ w.p. }  p_V,\\
b_V & \quad \text{ w.p. }  1-p_V,
\end{cases}
\end{equation}
with $0<a_U, a_V < 1 < b_U, b_V$. It follows that
\[
p_U = \frac{b_U-1}{b_U-a_U}, \sigma^2_U = (1-a_U)(b_U-1), c_3(U)=(a_U+b_U-2)\sigma_U^2,\]
and similar for $p_V, \sigma^2_V$, and $c_3(V)$.
The LSTs are given by
\begin{align}
\dE\big[e^{-\zeta V}\big]&=\frac{(b_V-1)e^{-\zeta a_V}+(1-a_V)e^{-\zeta b_V}}{b_V-a_V},\\
\dE\big[e^{-\zeta U/\rho}\big]&=\frac{(b_U-1)e^{-\zeta a_U/\rho}+(1-a_U)e^{-\zeta b_U/\rho}}{b_U-a_U}.
\end{align}

The interesting feature of this example is that, with $U$ and $V$ having a lattice distribution, Assumption~\ref{ass:delta} of Section \ref{sec5}, required for the Nearly-deterministic Gaussian regime, is \emph{not} (necessarily) satisfied. Nevertheless, we still get accurate approximations when applying our saddle point technique, as shown in Table~\ref{tbl:3-1f}, where we take $a_U=3/4$, $b_U=2$, $a_V=1/2$, $b_V=3/2$. For completeness, we also present the results for the classical Kingman-type heavy traffic regime and for the nearly-deterministic Kingman-type heavy traffic regime.

\noindent\textbf{Classical HT Kingman.} Again, we approximate the moments of $\alpha W$ by \eqref{eqn80}, with
\beq \label{6.9c}
\sigma_{\alpha}^2=(\sigma_V^2+\rho^{-2}\sigma_U^2)\,\rho=~\rho\left((1-a_V)(b_V-1)+\frac{(1-a_U)(b_U-1)}{\rho^2}\right),
\eq
and $\rho=1-\alpha$, where we take $\alpha=1/10,1/100,1/1000$. The results, for the first five moments of $W$, are presented in Table~\ref{tbl:1-1f}.

\begin{table}[t]
\[
\begin{array}{|l|rr|rr|rr|}
%\hline
%\mc{7}{c}{\text{Example 3  - Classical HT Kingman}~~ (\vart_U=5/2, m=4, \delta=1)}\\
\hline
& \mc{2}{c|}{\alpha=1/10} & \mc{2}{c|}{\alpha=1/100} & \mc{2}{c|}{\alpha=1/1000} \\
\hline
k & {\rm Exact} & {\rm Asymp} & {\rm Exact} & {\rm Asymp} & {\rm Exact} & {\rm Asymp} \\
\hline
 1 & 0.206 & 0.251 & 0.246 & 0.250 & 0.250 & 0.250 \\
 2 & 0.097 & 0.126 & 0.122 & 0.125 & 0.125 & 0.125 \\
 3 & 0.068 & 0.095 & 0.091 & 0.094 & 0.093 & 0.094 \\
 4 & 0.064 & 0.096 & 0.091 & 0.094 & 0.093 & 0.094 \\
 5 & 0.075 & 0.120 & 0.113 & 0.117 & 0.117 & 0.117 \\
\hline
\end{array}
\]
\caption{Example 4  - Classical HT Kingman: Comparison of exact result for $m_k(\alpha W)$ %, obtained via (\ref{6.7}) and numerical computation of the cumulants in (\ref{6.10}),
and asymptotic result $k!(\frac12\,\sigma_{\alpha}^2)^k$% with $\sigma_{\alpha}^2$ given in (\ref{6.9})
%, for the case that %$\vart_U=5/2$, $m=4$, $\delta = 1$
with $\alpha=1/10,1/100,1/1000$ and $k=1,2,3,4,5$.}
\label{tbl:1-1f}
\end{table}

\noindent\textbf{Nearly-deterministic HT Kingman.} With $\beta=1$ and $\rho_n=1-\beta/n$, we get
\beq \label{6.12b}
\gamma_n=(\sigma_V^2+\rho_n^{-2}\sigma_U^2)\,\rho_n/(2n(1-\rho_n))=\frac{1+\rho_n^2}{8n\rho_n(1-\rho_n)}~.
\eq
As in Example 2, we see that $\gamma_n$ in the Nearly-deterministic HT Kingman regime is equal to $\frac12\sigma_{\alpha}^2$, due to the fact that $\beta=1$, resulting in the same limiting values (shown in Table~\ref{tbl:2-1f}).

\begin{table}[t]
\[
\begin{array}{|l|rr|rr|rr|rr|rr|}
%\hline
%\mc{7}{c}{\text{Example 3  - Nearly-deterministic HT Kingman}~~ (\vart_U=5/2, m=4, \delta=1)}\\
\hline
& \mc{2}{c|}{n=10} & \mc{2}{c|}{n=100} & \mc{2}{c|}{n=1000} & \mc{2}{c|}{n=10000}  & \mc{2}{c|}{n=100000} \\
\hline
k & {\rm Exact} & {\rm Asymp} & {\rm Exact} & {\rm Asymp} & {\rm Exact} & {\rm Asymp} & {\rm Exact} & {\rm Asymp}  & {\rm Exact} & {\rm Asymp} \\
\hline
 1 & 0.134 & 0.251 & 0.211 & 0.250 & 0.237 & 0.250 & 0.246 & 0.250 & 0.249 & 0.250 \\
 2 & 0.065 & 0.126 & 0.105 & 0.125 & 0.119 & 0.125 & 0.123 & 0.125 & 0.124 & 0.125 \\
 3 & 0.046 & 0.095 & 0.078 & 0.094 & 0.089 & 0.094 & 0.092 & 0.094 & 0.093 & 0.094 \\
 4 & 0.043 & 0.096 & 0.078 & 0.094 & 0.089 & 0.094 & 0.092 & 0.094 & 0.093 & 0.094 \\
 5 & 0.051 & 0.120 & 0.097 & 0.117 & 0.111 & 0.117 & 0.115 & 0.117 & 0.117 & 0.117 \\
\hline
\end{array}
\]
\caption{Example 4  - Nearly-deterministic HT Kingman: Comparison of exact result for $m_k(W_n)$ %, obtained via (\ref{6.7}) and numerical computation of the cumulants in (\ref{6.10}),
and asymptotic result $k!\gamma_n^k$% with $\sigma_{\alpha}^2$ given in (\ref{6.9})
%, for the case that $\vart_U=5/2$, $m=4$, $\delta = 1$
 with $n=10, 100, 1000, 10000, 100000$ and $k=1, 2, 3, 4, 5$.}
\label{tbl:2-1f}
\end{table}

\noindent\textbf{Nearly-deterministic HT Gaussian.}
To use the saddle point method, we first compute% $h(\zeta)$:
\begin{align}
h(\zeta)%=\log \psi(-\zeta)
&= \log\left[\frac{(b_V-1)e^{-\zeta a_V}+(1-a_V)e^{-\zeta b_V}}{b_V-a_V}\right]+
\log\left[\frac{(b_U-1)e^{-\zeta a_U/\rho}+(1-a_U)e^{-\zeta b_U/\rho}}{b_U-a_U}\right].
\end{align}
As in the previous example, there is no closed-from expression for $\zeta_{sp}$, but with all the derivatives of $h(\zeta)$ being analytic, numerical methods can easily find the saddle point of $h(\zeta)$.   It is noteworthy that, for our chosen parameter values, the function $h(\zeta_{sp}+iy)$, $y\in\dR$, that occurs as
\begin{equation}
\log (1-\psi^n(-\zeta))=\log(1-e^{nh(\zeta)}), \qquad\zeta=\zeta_{sp}+iy
\end{equation}
in \eqref{5.8} and \eqref{eqn85} is \emph{periodic} in $y$, with period $8k\pi$ as $\rho=1-1/k$, $k=2,3,\dots$. Hence, the maximum value $|h(\zeta_{sp})|$ of $|h(\zeta_{sp}+iy)|$ is repeated infinitely many times (far) away from the negative real axis.
Still, our numerical methods encountered no issues computing the exact values and the two approximations for the first five moments of $\frac{\sqrt{n}}{\sigma_n}W_n$. In this example, with $c_3(V)=0$ and $c_3(U)=3/16$, the refined approximation gives significantly better results than the standard approximation  for $k=1$ and $k=2$. For $k\geq 3$, the improvement less significant. The results are shown in Table~\ref{tbl:3-1f}.

%In this example, with $c_3(V)=c_3(U)=0$, the refined approximation gives significantly better results than the standard approximation.
\begin{table}[t]
\[
\begin{array}{|l|rrr|rrr|rrr|}
%\hline
%\mc{10}{c}{\text{Example 3  - Nearly-deterministic HT Gaussian}~~ (\vart_U=5/2,  m=4, \delta=1)}\\
\hline
& \mc{3}{c|}{n=10} & \mc{3}{c|}{n=100} & \mc{3}{c|}{n=1000} \\
\hline
k & {\rm Exact} & {\rm Asymp\ 1} & {\rm Asymp\ 2} & {\rm Exact} & {\rm Asymp\ 1} & {\rm Asymp\ 2} & {\rm Exact} & {\rm Asymp\ 1} & {\rm Asymp\ 2} \\
\hline
 1 & 0.0186 & 0.0270 & 0.0190 & 0.0350 & 0.0362 & 0.0349 & 0.0406 & 0.0409 & 0.0406 \\
 2 & 0.0149 & 0.0218 & 0.0149 & 0.0295 & 0.0307 & 0.0294 & 0.0352 & 0.0355 & 0.0352 \\
 3 & 0.0157 & 0.0244 & 0.0147 & 0.0344 & 0.0361 & 0.0343 & 0.0419 & 0.0428 & 0.0423 \\
 4 & 0.0202 & 0.0345 & 0.0147 & 0.0505 & 0.0538 & 0.0502 & 0.0692 & 0.0652 & 0.0642 \\
 5 & 0.0298 & 0.0584 & 0.0077 & 0.0890 & 0.0961 & 0.0878 & 0.1182 & 0.1196 & 0.1173 \\
\hline
\end{array}
\]
\caption{Example 4  - Nearly-deterministic HT Gaussian: Comparison of exact result for $m_k(\frac{\sqrt{n}}{\sigma_n}\, W_n)$ and the asymptotic results% for the case that $\vart_U=5/2$, $m=4$, $\delta = 1$ with $\beta=1$,
with $n=10, 100, 1000$ and $k=1, 2, 3, 4, 5$.  The two entries in the Asymp-columns give, for a particular $k$, the asymptotic result from (\ref{eqn86}) and (\ref{eqn89}) in that order.}
\label{tbl:3-1f}
\end{table}

\subsection{Example 5: $U$ Inverse Gaussian distributed and $V$ Gamma distributed}

We take $U\sim IG(\mu, \lambda)$, inverse Gaussian with parameters $\mu=1$ and $\lambda=\frac{10}{8}\sqrt{\frac32}$, and $V\sim \text{Gamma}(k_V, \vart_V)$ with $\vart_V=1/k_V=8/10$. The probability density function of $U$ is given by
\begin{equation}\label{eqn:invgauss}
\sqrt{\frac{\lambda}{2\pi x^3}}\,\exp\left[-\frac{\lambda(x-\mu)^2}{2\mu^2x}\right], \qquad x\geq 0; \lambda,\mu >0.
\end{equation}
For this particular combination of distribution functions the saddle point $\zeta_{sp}$ has an explicit form, viz.
\begin{align}
h(\zeta) &= -\frac{1}{\vart_V}\log(1+\vart_V \zeta)+\lambda\big(1-\sqrt{1-{2\zeta}/{\lambda\rho}}\big),\label{ex4eq1}\\
\zeta_{sp}&= \frac{1}{\vart_V}\left(\sqrt{(\gamma+1)^2-(1-\rho^2)}-\gamma-1\right),\qquad \text{ with }\gamma=\frac{\rho}{\lambda\vart_V}.\label{ex4eq2}
\end{align}

The parameters have been chosen such that $\vart_V \lambda=\sqrt{3/2}$, which results in the third cumulants of $U$ and $V$ being equal:
\[
 c_3(U)=\frac{3}{\lambda^2}=\frac{32}{25}, \quad c_3(V)=2\vart^2=\frac{32}{25}.
\]
The purpose of this example is to show that the accuracy of the two approximations in estimating the moments of the scaled waiting time distribution is comparable in the nearly deterministic Gaussian regime. Recall the third case in Example~1, where $h'''(\zeta_{sp})$ was equal to zero, causing the two approximations to be \emph{exactly} the same. Now we have two different distributions for $U$ and $V$, with $h'''(\zeta_{sp})\rightarrow 0$ as $n\to\infty$. As a consequence, the two approximations are only the same in the limit.

For completeness, we also provide the tables for the classical HT Kingman case and the nearly-deterministic HT Kingman case.

\noindent\textbf{Classical HT Kingman.}
We approximate the moments of $\alpha W$ by $k! (\frac12 \sigma_\alpha^2)^k$, with
\begin{equation}
\sigma_{\alpha}^2=(\sigma_V^2+\rho^{-2}\sigma_U^2)\,\rho=\left(k_V\vart_V^2+\rho^{-2}
\frac{\mu^3}{\lambda}\right)\,\rho=(\rho\vart_V+1/\lambda\rho).
\end{equation}
The results are shown in Table \ref{tbl:1-1d}.
\begin{table}[t]
\[
\begin{array}{|l|rr|rr|rr|}
%\hline
%\mc{7}{c}{\text{Example 4  - Classical HT Kingman}~~ (\vart_V=1/10, \mu=1, \lambda=1/10)}\\
\hline
& \mc{2}{c|}{\alpha=1/10} & \mc{2}{c|}{\alpha=1/100} & \mc{2}{c|}{\alpha=1/1000} \\
\hline
k & {\rm Exact} & {\rm Asymp} & {\rm Exact} & {\rm Asymp} & {\rm Exact} & {\rm Asymp} \\
\hline
 1 & 0.634 & 0.723 & 0.717 & 0.726 & 0.726 & 0.727 \\
 2 & 0.912 & 1.045 & 1.041 & 1.054 & 1.054 & 1.056 \\
 3 & 1.969 & 2.267 & 2.268 & 2.295 & 2.298 & 2.301 \\
 4 & 5.666 & 6.554 & 6.584 & 6.664 & 6.679 & 6.687 \\
 5 & 20.381 & 23.688 & 23.894 & 24.186 & 24.262 & 24.290 \\
\hline
\end{array}
\]
\caption{Example 5  - Classical HT Kingman: Comparison of exact result for $m_k(\alpha W)$ %, obtained via (\ref{6.7}) and numerical computation of the cumulants in (\ref{6.10}),
and asymptotic result $k!(\frac12\,\sigma_{\alpha}^2)^k$ % with $\sigma_{\alpha}^2$ given in (\ref{6.9})
%, for the case that $\vart_V=1/10$, $\mu=1$, $\lambda=1/10$
 with $\alpha=1/10,1/100,1/1000$ and $k=1,2,3,4,5$.}
\label{tbl:1-1d}
\end{table}

\noindent\textbf{Nearly-deterministic HT Kingman.}
In the numerical results, we take $\beta=1$ again and approximate $\dE[W_n^k]$ by $k!\gamma_n^k$, with
\beq
\gamma_n=(\sigma_V^2+\rho_n^{-2}\sigma_U^2)\,\rho_n/(2n(1-\rho_n))=(\rho_n\vart_V+1/{\lambda\rho_n})\,/(2n(1-\rho_n))~.
\eq
The results are shown in Table \ref{tbl:2-1d}.

\begin{table}[t]
\[
\begin{array}{|l|rr|rr|rr|rr|rr|}
%\hline
%\mc{7}{c}{\text{Example 4  - Nearly-deterministic HT Kingman}~~ (\vart_V=1/10, \mu=1, \lambda=1/10)}\\
\hline
& \mc{2}{c|}{n=10} & \mc{2}{c|}{n=100} & \mc{2}{c|}{n=1000} & \mc{2}{c|}{n=10000} & \mc{2}{c|}{n=100000} \\
\hline
k & {\rm Exact} & {\rm Asymp} & {\rm Exact} & {\rm Asymp} & {\rm Exact} & {\rm Asymp}& {\rm Exact} & {\rm Asymp} & {\rm Exact} & {\rm Asymp} \\
\hline
 1 & 0.509 & 0.723 & 0.658 & 0.726 & 0.705 & 0.727 & 0.720 & 0.727 & 0.724 & 0.727 \\
 2 & 0.740 & 1.045 & 0.956 & 1.054 & 1.024 & 1.056 & 1.046 & 1.056 & 1.053 & 1.056 \\
 3 & 1.601 & 2.267 & 2.082 & 2.295 & 2.232 & 2.301 & 2.279 & 2.302 & 2.295 & 2.302 \\
 4 & 4.608 & 6.554 & 6.044 & 6.664 & 6.485 & 6.687 & 6.625 & 6.689 & 6.669 & 6.689 \\
 5 & 16.578 & 23.688 & 21.935 & 24.186 & 23.559 & 24.290 & 24.068 & 24.301 & 24.228 & 24.303 \\
\hline
\end{array}
\]
\caption{Example 5  - Nearly-deterministic HT Kingman: Comparison of exact result for $m_k(W_n)$ %, obtained via (\ref{6.7}) and numerical computation of the cumulants in (\ref{6.10}),
and asymptotic result $k!\gamma_n^k$ % with $\sigma_{\alpha}^2$ given in (\ref{6.9})
%, for the case that $\vart_V=1/10$, $\mu=1$, $\lambda=1/10$
 with $n=10, 100, 1000, 10000, 100000$ and $k=1, 2, 3, 4, 5$.}
\label{tbl:2-1d}
\end{table}

\noindent\textbf{Nearly-deterministic HT Gaussian.}
We invoke the machinery of the dedicated saddle point method with $h(\zeta)$ and $\zeta_{sp}$ as given in \eqref{ex4eq1} and \eqref{ex4eq2}, respectively. The results, shown in Table \ref{tbl:3-1d}, clearly indicate that the difference between the two approximations quickly vanishes as $n$ increases. It is striking to see that the approximations seem to give extremely accurate results in this case. In particular, it seems that the standard approximation slightly outperforms the enhanced version, but both are so close that it hardly makes a difference. The effect of taking $c_3(U)=c_3(V)$ can be observed very clearly in this example, in the sense that now the standard approximation Asymp 1 slightly outperforms the refined approximation Asymp 2.
\begin{table}[t]
\[
\begin{array}{|l|rrr|rrr|rrr|}
%\hline
%\mc{10}{c}{\text{Example 5  - Nearly-deterministic HT Gaussian}~~ (\vart_U=5/2,  m=4, \delta=1)}\\
\hline
& \mc{3}{c|}{n=10} & \mc{3}{c|}{n=100} & \mc{3}{c|}{n=1000} \\
\hline
k & {\rm Exact} & {\rm Asymp\ 1} & {\rm Asymp\ 2} & {\rm Exact} & {\rm Asymp\ 1} & {\rm Asymp\ 2} & {\rm Exact} & {\rm Asymp\ 1} & {\rm Asymp\ 2} \\
\hline
% 1 & 0.119 & 0.121 & 0.122 & 0.174 & 0.174 & 0.174 & 0.191 & 0.191 & 0.191 \\
% 2 & 0.132 & 0.136 & 0.138 & 0.220 & 0.220 & 0.220 & 0.250 & 0.250 & 0.250 \\
% 3 & 0.210 & 0.215 & 0.220 & 0.393 & 0.394 & 0.395 & 0.463 & 0.463 & 0.463 \\
% 4 & 0.431 & 0.436 & 0.453 & 0.913 & 0.914 & 0.917 & 1.115 & 1.115 & 1.116 \\
% 5 & 1.087 & 1.086 & 1.140 & 2.618 & 2.618 & 2.629 & 3.320 & 3.320 & 3.321 \\
%
  1 & 0.1185 & 0.1212 & 0.1222 & 0.1740 & 0.1744 & 0.1745 & 0.1913 & 0.1913 & 0.1913 \\
 2 & 0.1323 & 0.1359 & 0.1381 & 0.2196 & 0.2202 & 0.2204 & 0.2497 & 0.2497 & 0.2497 \\
 3 & 0.2099 & 0.2147 & 0.2203 & 0.3932 & 0.3941 & 0.3948 & 0.4630 & 0.4631 & 0.4632 \\
 4 & 0.4306 & 0.4360 & 0.4525 & 0.9131 & 0.9144 & 0.9170 & 1.1151 & 1.1153 & 1.1156 \\
 5 & 1.0868 & 1.0863 & 1.1404 & 2.6176 & 2.6184 & 2.6287 & 3.3196 & 3.3197 & 3.3209 \\
\hline
\end{array}
\]
\caption{Example 5  - Nearly-deterministic HT Gaussian: Comparison of exact result for $m_k(\frac{\sqrt{n}}{\sigma_n}\, W_n)$ and the asymptotic results %for the case that $\vart_V=1/10$, $\mu=1$, $\lambda=1/10$
with $\beta=1$, $n=10, 100, 1000$ and $k=1, 2, 3, 4, 5$.  The two entries in the Asymp-columns give, for a particular $k$, the asymptotic result from (\ref{eqn86}) and (\ref{eqn89}) in that order.}
\label{tbl:3-1d}
\end{table}

\subsection{Example 6: distributional results}

\begin{figure}[ht]
\parbox{0.48\linewidth}{%
\centering
\includegraphics[width=\linewidth]{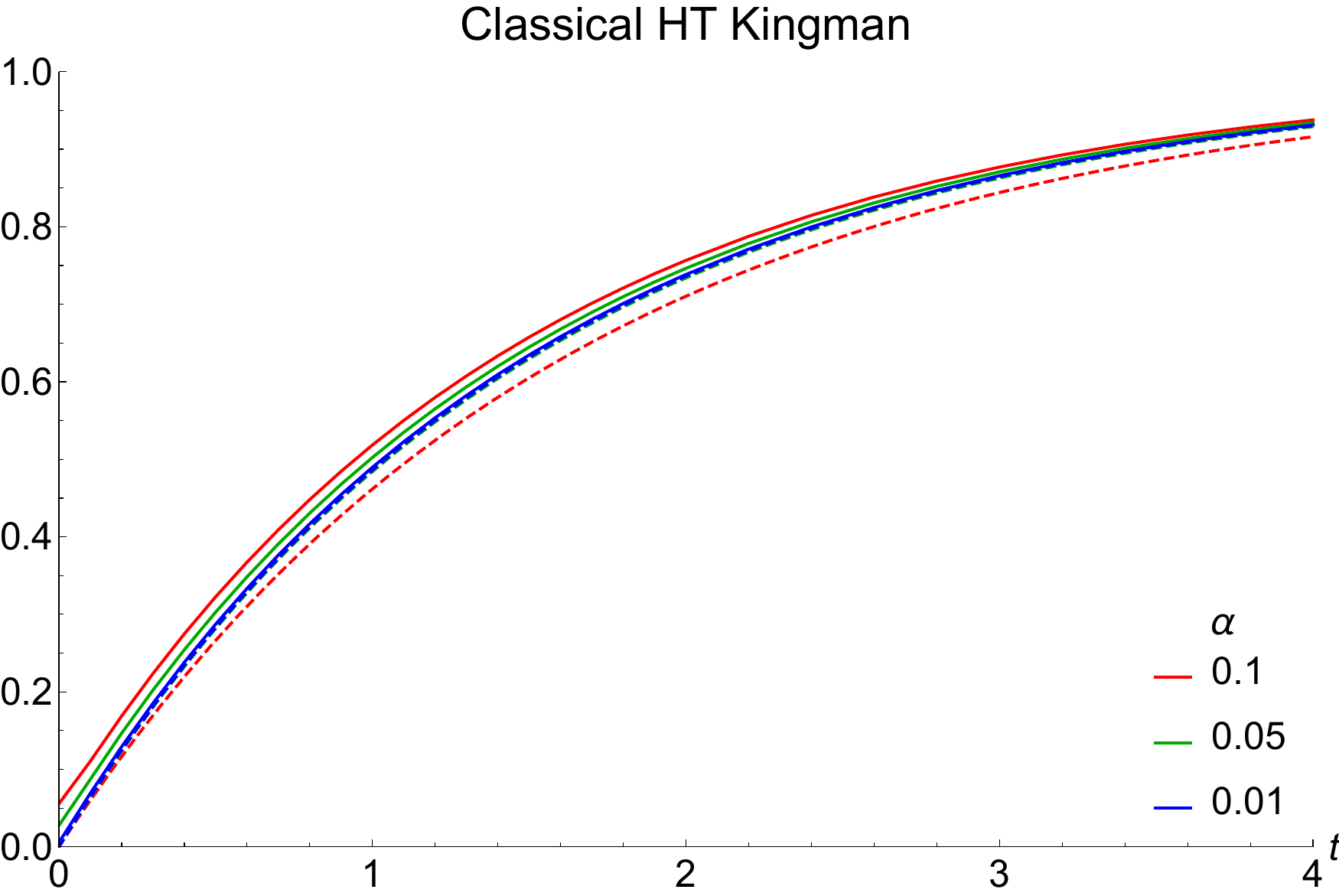}\\[-1ex]
\small (a) $P(\alpha W < t)$ %Classical HT Kingman
}
\hfill
\parbox{0.48\linewidth}{%
\centering
\includegraphics[width=\linewidth]{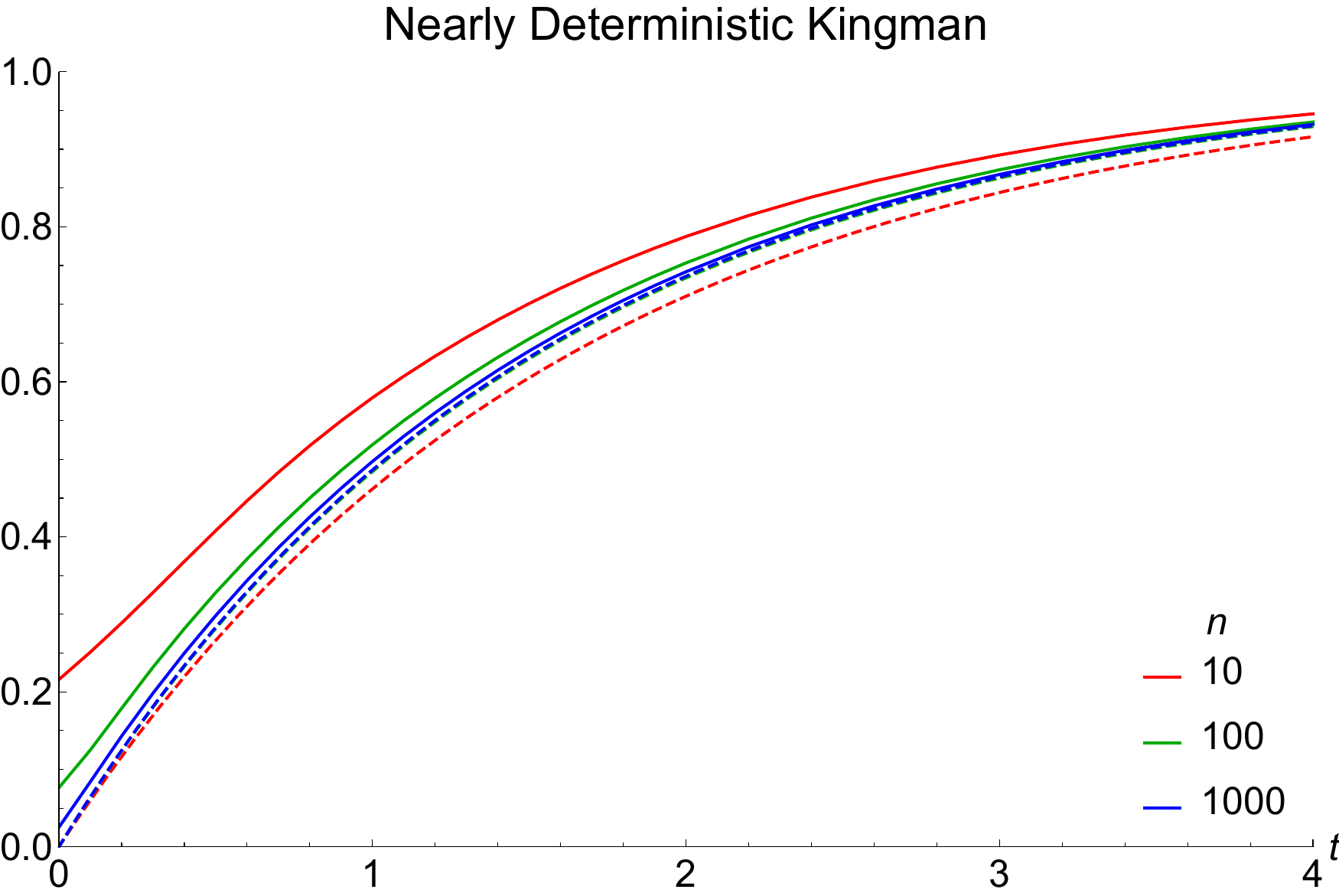}\\[-1ex]
\small (b) $P(W_n < t)$ % Nearly deterministic Kingman
}\\[1ex]
\parbox{0.48\linewidth}{%
\centering
\includegraphics[width=\linewidth]{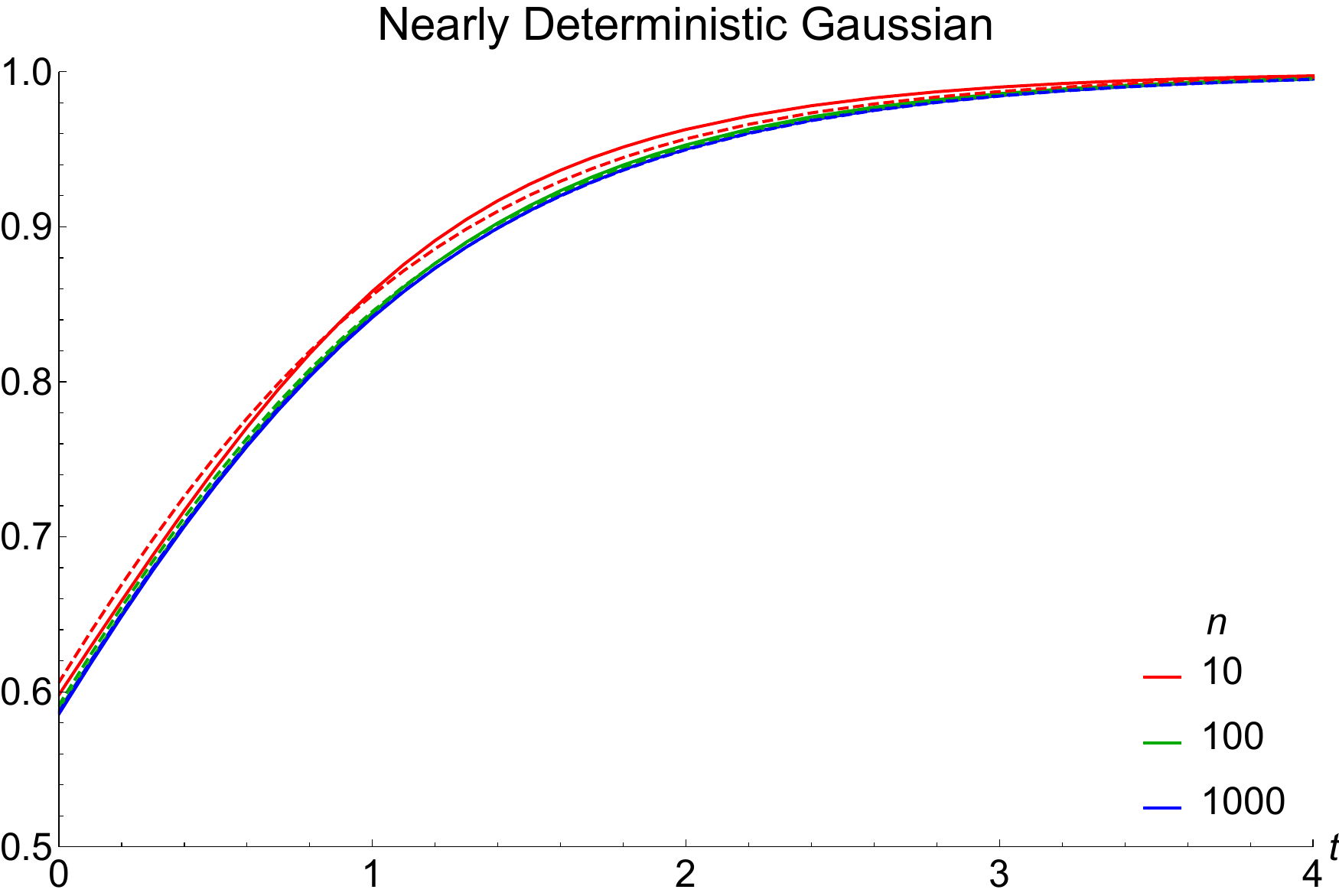}\\[-1ex]
\small (c) $P(\frac{\sqrt{n}}{\sigma_n} W_n < t)$ % Nearly deterministic Gaussian\\(standard approximation)
}
\hfill
\parbox{0.48\linewidth}{%
\centering
\includegraphics[width=\linewidth]{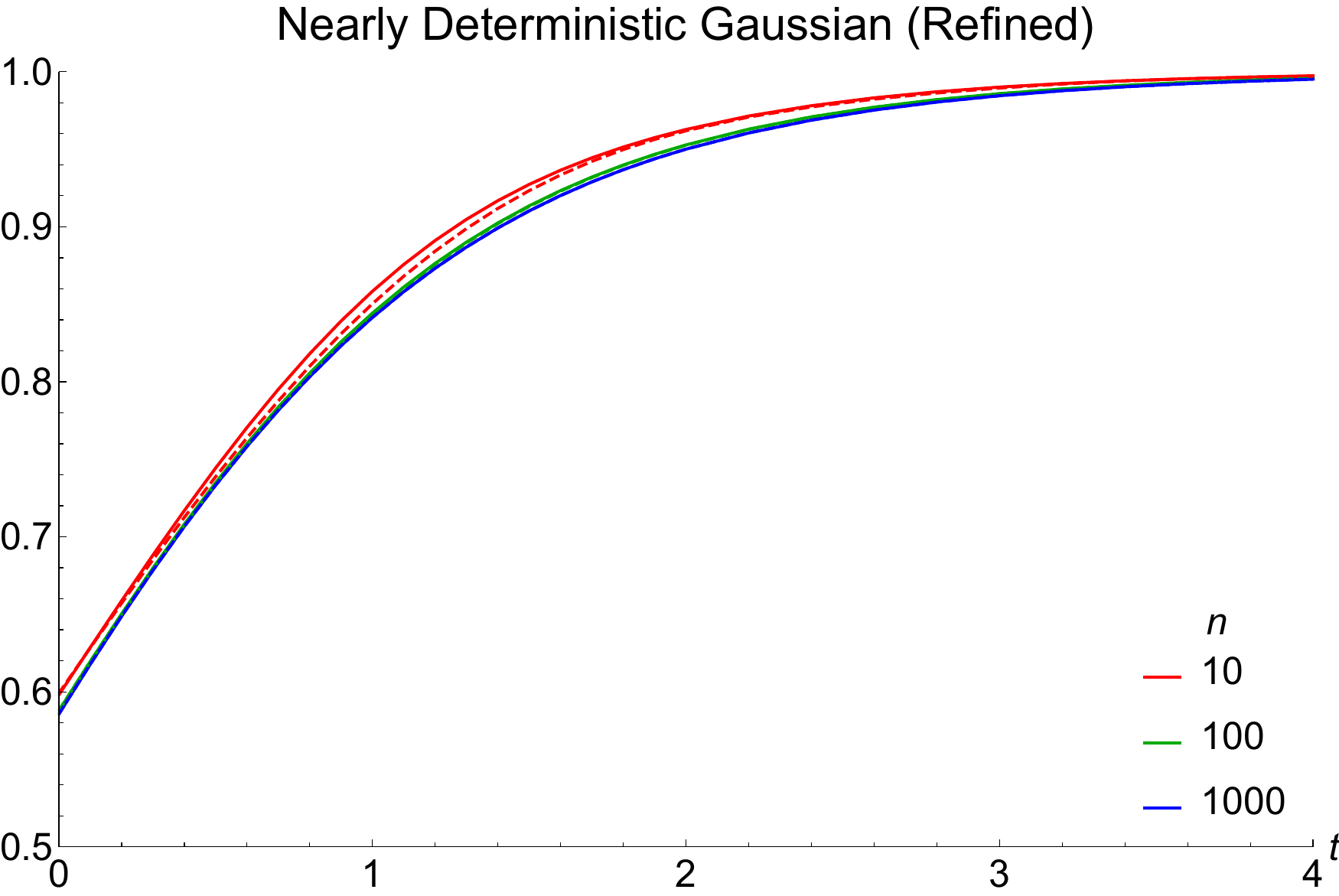}\\[-1ex]
\small (d) $P(\frac{\sqrt{n}}{\sigma_n} W_n < t)$ %Nearly deterministic Gaussian\\(refined approximation)
}
\caption{Example 6: CDFs of the scaled waiting-time distributions.}
\label{fig:cdfs}
\end{figure}

While our approach to obtain the convergence result is, in principle, meant for the moments of the (scaled) waiting times, we experienced no issues when trying to invert the Laplace-Stieltjes transforms to obtain the cumulative distribution functions of the waiting times. We used Abate and Whitt's numerical LST inversion algorithm \cite[EULER algorithm]{abatewhitt92}, which evaluates the LST only at complex points with a nonnegative real part. In Figures~\ref{fig:cdfs}(a)--(d) we have plotted the CDFs of the steady-state scaled waiting-time distribution for the first model discussed in Example~1 ($U$ and $V$ both Gamma distributed, with $\vart_U=5/2$ and $\vart_V=1/2$). The solid lines represent the exact CDFs, obtained by numerical inversion of \eqref{75}, appropriately scaled, with $\psi(-z)$ replaced by $\phi(-z)$ for the nearly deterministic regimes.

In Fig.~\ref{fig:cdfs}(a), results are shown for the classical Kingman-type heavy traffic, for $\alpha=0.1$, $\alpha=0.05$, and $\alpha=0.01$, respectively. The dashed lines are the corresponding approximations based on the limiting exponential distribution as given in Theorem~\ref{thm1}. Figure~\ref{fig:cdfs}(b) gives the exact results for the nearly deterministic Kingman case, with $n=10, 100, 1000$. The approximations are the CDFs of the exponential distributions given in Theorem~\ref{thm2}. Figures~\ref{fig:cdfs}(c)--(d) are similar, but for the nearly-deterministic Gaussian regimes. For these plots, we inverted the LST of $M_\beta$ given in Equation~\eqref{5.12} to obtain its CDF. The standard approximation based on Theorem~\ref{thm3} is shown in Fig.~\ref{fig:cdfs}(c), whereas the refined approximation based on Theorem~\ref{thm5} is shown in Fig.~\ref{fig:cdfs}(d).

It is clear that the approximations for the classical and nearly deterministic Kingman cases, which are based on the exponential distribution, perform worse for smaller values of $\rho$. This is most prominently visible when the probability of having a waiting time of zero is substantial. Compare, for example, Figures~\ref{fig:cdfs}(a) and \ref{fig:cdfs}(b). In the first case,
corresponding to the classical Kingman HT regime, when $\alpha=0.1$ we find $P(W=0) = 0.056$. In the second case, in the nearly-deterministic Kingman HT regime with $n=10$, we find $P(W=0) = 0.216$. Obviously, this probability is zero in all our approximations, meaning that the approximation performs much worse in the classical HT regime than in the nearly-deterministic HT regime, despite the fact that in both cases the load of the system is $\rho=0.9$. In contrast, the approximations for the Gaussian nearly-deterministic HT regime are extremely accurate for $P(W=0)$, even for lower values of $n$.

\subsection{Example 7: The $M/D/s$ queue}

In this example, we show that the techniques can also be used to analyze the $M/D/s$ queue. This multiserver queue fits within our framework due to the equivalence between the $M/D/s$ queue and the \emph{single}-server queue with Erlang$(s)$ interarrival times and deterministic service times,  with respect to the waiting times (cf. \cite{mds}).

\noindent\textbf{Classical HT Kingman.} For this model, we have \eqref{eqn80}
\beq \label{6.9b}
\sigma_{\alpha}^2=(\sigma_V^2+\rho^{-2}\sigma_U^2)\,\rho=1/(s\rho)~.
\eq
With $\rho=1-\alpha$ and $\alpha=1/10,1/100,1/1000$, we see that the approximation for the first five moments of $W$, presented in Table~\ref{tbl:1-1x}, is very accurate for values of $\alpha$ less than $1/100$, which corresponds to $\rho>0.99$.

\begin{table}[t]
\[
\begin{array}{|l|rr|rr|rr|}
%\hline
%\mc{7}{c}{\text{Example 2  - Classical HT Kingman}~~ (\vart_U=5/2, m=4, \delta=1)}\\
\hline
& \mc{2}{c|}{\alpha=1/10} & \mc{2}{c|}{\alpha=1/100} & \mc{2}{c|}{\alpha=1/1000} \\
\hline
k & {\rm Exact} & {\rm Asymp} & {\rm Exact} & {\rm Asymp} & {\rm Exact} & {\rm Asymp} \\
\hline
 1 & 0.078 & 0.111 & 0.098 & 0.101 & 0.100 & 0.100 \\
 2 & 0.015 & 0.025 & 0.019 & 0.020 & 0.020 & 0.020 \\
 3 & 0.004 & 0.008 & 0.006 & 0.006 & 0.006 & 0.006 \\
 4 & 0.002 & 0.004 & 0.002 & 0.002 & 0.002 & 0.002 \\
 5 & 0.001 & 0.002 & 0.001 & 0.001 & 0.001 & 0.001 \\
\hline
\end{array}
\]
\caption{Example 2  - Classical HT Kingman: Comparison of exact result for $m_k(\alpha W)$ %, obtained via (\ref{6.7}) and numerical computation of the cumulants in (\ref{6.10}),
and asymptotic result $k!(\frac12\,\sigma_{\alpha}^2)^k$% with $\sigma_{\alpha}^2$ given in (\ref{6.9})
%, for the case that $\vart_U=5/2$, $m=4$, $\delta = 1$
 with $\alpha=1/10,1/100,1/1000$ and $k=1,2,3,4,5$.}
\label{tbl:1-1x}
\end{table}

\noindent\textbf{Nearly-deterministic HT Kingman.} With $\beta=1$ and $\rho_n=1-\beta/n$, we get
\beq \label{6.12b}
\gamma_n=(\sigma_V^2+\rho_n^{-2}\sigma_U^2)\,\rho_n/(2n(1-\rho_n))=
\frac{1}{2sn\rho_n(1-\rho_n)}.
~.
\eq
In Table \ref{tbl:2-1x}, we show results for for $n=10, 100, 1000, 10000, 100000$, which corresponds to $\rho_n=0.9, 0.99$, $0.999$, $0.9999$, $0.99999$, respectively.

\begin{table}[t]
\[
\begin{array}{|l|rr|rr|rr|rr|rr|}
%\hline
%\mc{7}{c}{\text{Example 2  - Nearly-deterministic HT Kingman}~~ (\vart_U=1/2, m=4, \delta=1)}\\
\hline
& \mc{2}{c|}{n=10} & \mc{2}{c|}{n=100} & \mc{2}{c|}{n=1000}& \mc{2}{c|}{n=10000} & \mc{2}{c|}{n=100000} \\
\hline
k & {\rm Exact} & {\rm Asymp} & {\rm Exact} & {\rm Asymp} & {\rm Exact} & {\rm Asymp} & {\rm Exact} & {\rm Asymp} & {\rm Exact} & {\rm Asymp} \\
\hline
 1 & 0.040 & 0.111 & 0.077 & 0.101 & 0.092 & 0.100 & 0.097 & 0.100 & 0.099 & 0.100 \\
 2 & 0.008 & 0.025 & 0.015 & 0.020 & 0.018 & 0.020 & 0.019 & 0.020 & 0.020 & 0.020 \\
 3 & 0.002 & 0.008 & 0.005 & 0.006 & 0.006 & 0.006 & 0.006 & 0.006 & 0.006 & 0.006 \\
 4 & 0.001 & 0.004 & 0.002 & 0.002 & 0.002 & 0.002 & 0.002 & 0.002 & 0.002 & 0.002 \\
 5 & 0.000 & 0.002 & 0.001 & 0.001 & 0.001 & 0.001 & 0.001 & 0.001 & 0.001 & 0.001 \\
\hline
\end{array}
\]
\caption{Example 2  - Nearly-deterministic HT Kingman: Comparison of exact result for $m_k(\alpha W)$ %, obtained via (\ref{6.7}) and numerical computation of the cumulants in (\ref{6.10}),
and asymptotic result $k!(\frac12\,\sigma_{\alpha}^2)^k$% with $\sigma_{\alpha}^2$ given in (\ref{6.9})
%, for the case that $\vart_U=5/2$, $m=4$, $\delta = 1$
 with $n=10,100,1000,10000,100000$ and $k=1,2,3,4,5$.}
\label{tbl:2-1x}
\end{table}

\noindent\textbf{Nearly-deterministic HT Gaussian.}
It is readily seen that
\begin{equation}
h(\zeta)= -s \log(1-\zeta/(s\rho_n))-\zeta.
\end{equation}

The exact and approximate values for the first five moments of the scaled waiting times are shown in Table \ref{tbl:3-1x}.

\begin{table}[t]
\[
\begin{array}{|l|rrr|rrr|rrr|}
%\hline
%\mc{10}{c}{\text{Example 2  - Nearly-deterministic HT Gaussian}~~ (\vart_U=1/2,  m=4, \delta=1)}\\
\hline
& \mc{3}{c|}{n=10} & \mc{3}{c|}{n=100} & \mc{3}{c|}{n=1000} \\
\hline
k & {\rm Exact} & {\rm Asymp\ 1} & {\rm Asymp\ 2} & {\rm Exact} & {\rm Asymp\ 1} & {\rm Asymp\ 2} & {\rm Exact} & {\rm Asymp\ 1} & {\rm Asymp\ 2} \\
\hline
 1 & 0.0024 & 0.0046 & 0.0026 & 0.0039 & 0.0046 & 0.0039 & 0.0043 & 0.0046 & 0.0043 \\
 2 & 0.0015 & 0.0030 & 0.0017 & 0.0025 & 0.0030 & 0.0025 & 0.0028 & 0.0030 & 0.0028 \\
 3 & 0.0014 & 0.0027 & 0.0014 & 0.0022 & 0.0027 & 0.0023 & 0.0025 & 0.0027 & 0.0025 \\
 4 & 0.0015 & 0.0031 & 0.0014 & 0.0025 & 0.0031 & 0.0025 & 0.0029 & 0.0031 & 0.0029 \\
 5 & 0.0019 & 0.0042 & 0.0015 & 0.0034 & 0.0042 & 0.0033 & 0.0000 & 0.0042 & 0.0039 \\
\hline
\end{array}
\]
\caption{Example 2  - Nearly-deterministic HT Gaussian: Comparison of exact result for $m_k(\frac{\sqrt{n}}{\sigma_n}\, W_n)$ and the asymptotic results for %the case% that $\vart_U=5/2$, $m=4$, $\delta = \frac12$
% with $\beta=1$,
 $n=10,100,1000$ and $k=1,2,3,4,5$.  The two entries in the Asymp-columns give, for a particular $k$, the asymptotic result from (\ref{eqn86}) and (\ref{eqn89}) in that order.}
\label{tbl:3-1x}
\end{table}

\end{document}